\documentclass[11pt]{amsart}
\usepackage{amssymb,amsmath}
\usepackage{multicol}
\usepackage{amsfonts}
\usepackage{dsfont}
  \usepackage{paralist}
\usepackage{enumerate}
\usepackage[latin1]{inputenc}
\usepackage{color}
\usepackage[a4paper]{geometry}
\usepackage{graphicx}
\DeclareMathOperator{\Jac}{Jac}
\DeclareMathOperator{\dist}{dist}

\newtheorem{Prop}{Main Proposition}

\newtheorem{defi}{Definition}
\newtheorem{lem}{Lemma}
\newtheorem{Lem}{Reparametrization Lemma}

\newtheorem{prop}{Proposition}
\newtheorem{theo}{Theorem}

\newenvironment{demo}
{\medbreak\noindent{\sc Proof :}} {\hfill $\square$ \medbreak}

\newenvironment{demof}[1]
{\medbreak\noindent{\sc Proof of {#1} :}}
{\hfill$\square$\medbreak}

\begin{document}

\newcommand\N{{\mathbb N}}
\newcommand\R{{\mathbb R}}

\title{\textbf{Symbolic extensions in intermediate smoothness on surfaces}\\
\vspace{1cm}
\textbf{Extensions symboliques  en r\'egularit\'e interm\'ediaire sur les surfaces}}

\author{David Burguet, \\
{\it CMLA-ENS Cachan}\\
{\it 61 avenue du président Wilson}\\ {\it 94235 Cachan Cedex France}
}

\date{}
\maketitle

\pagestyle{myheadings} \markboth{\normalsize\sc David
Burguet}{\normalsize\sc Symbolic extensions in intermediate smoothness on surfaces}

\noindent \textbf{Abstract :} We prove that $\mathcal{C}^r$  maps with $r>1$ on a compact surface have symbolic extensions, i.e. topological extensions which are subshifts over a finite alphabet.  More precisely we give a sharp upper bound on the so-called symbolic extension entropy, which is the infimum of the topological entropies of all the symbolic extensions.  This answers positively a conjecture of S.Newhouse and T.Downarowicz in dimension two and improves a previous result of the author \cite{burinv}.\\

\noindent \textbf{R\'esum\'e :} Nous montrons que toute dynamique de classe  $\mathcal{C}^r$ avec $r>1$ sur une surface compacte admet une extension symbolique, i.e. une extension topologique qui est un sous-décalage à alphabet fini. Nous donnons plus précisément une borne (optimale) sur l'infimum de l'entropie topologique de toutes les extensions symboliques. Ceci répond positivement à une conjecture de S.Newhouse and T.Downarowicz en dimension deux et améliore un résultat précédent de l'auteur \cite{burinv}.

\section{Introduction}

By a dynamical system $(X,T)$ we mean a continuous map $T$ on a compact metrizable space $X$. One well studied class of dynamical systems are the symbolic ones, i.e. closed subset $Y$ of $\mathcal{A}^{\mathbb{Z}}$, with a finite alphabet $\mathcal{A}$, endowed with the shift $S$. Such a pair $(Y,S)$ is also called a subshift. Given a dynamical system $(X,T)$ one wonders if there exists a symbolic extension $(Y,S)$ of $(X,T)$, i.e. a subshift $(Y,S)$ along a continuous map $\pi:Y\rightarrow X$ such that $\pi\circ S=T\circ \pi$. We first observe that dynamical systems with symbolic extensions have necessarily  finite topological entropy. When a dynamical system has symbolic extensions  we are interested in minimizing their entropy. The topological symbolic extension entropy $h_{sex}(T)=\inf \{h_{top}(Y,S)$ :  $(Y,S)$  is a symbolic extension of $(X,T)\}$ estimates how the dynamical system $(X,T)$ differs from a symbolic extension from the point of view of entropy. The problem of the existence of symbolic extensions  leads to a deep theory of entropy which was developed mainly by M.Boyle and T.Downarowicz, who related the existence of symbolic extensions and their entropy with the convergence of the entropy of $(X,T)$ computed at finer and finer scales \cite{BD}.

By using a result of J.Buzzi \cite{Buz} involving Yomdin theory,  M.Boyle, D.Fiebig, U.Fiebig \cite{bff}  proved that $\mathcal{C}^{\infty}$ maps on a compact manifold  admit principal symbolic extensions, i.e. symbolic extensions which preserve the entropy of invariant measures \cite{bff}.  On the other hand $\mathcal{C}^1$ maps without symbolic extensions  have been built in several works \cite{ND}, \cite{Asa}, \cite{bur}. In the present paper we consider dynamical systems of intermediate smoothness, i.e. $\mathcal{C}^r$ maps $T$ on a compact manifold with $1<r<+\infty$ (we mean that $T$ admits a derivative or order $\left\lceil r-1\right\rceil$ which is $r-\left\lceil r-1\right\rceil$-Hölder). T.Downarowicz and  A.Maass  have  recently proved  that $\mathcal{C}^r$ maps of the interval $f:[0,1]\rightarrow [0,1]$ with $1<r<+\infty$ have symbolic extensions \cite{Dow}. More precisely they showed that  $h_{sex}(f)\leq h_{top}(f)+\frac{\log\|f'\|_{\infty}}{r-1}$. The author built explicit examples  \cite{bur}  proving this upper bound is sharp. Similar $\mathcal{C}^r$ examples with large symbolic extension entropy have been previously built by T.Downarowicz and S.Newhouse for diffeomorphisms in higher dimension \cite{ND}. The results of T.Downarowicz and  A.Maass have been extended by the author in any dimension to nonuniformly entropy expanding maps (i.e. $\mathcal{C}^1$ maps whose ergodic measures with positive entropy have nonnegative Lyapunov exponents) of class $\mathcal{C}^r$ with $1<r<+\infty$ \cite{superbur}. More recently the author also proved the existence of symbolic extensions for $\mathcal{C}^2$ surface local  diffeomorphisms \cite{burinv}.  T.Downarowicz and S.Newhouse have conjectured in \cite{ND} that $\mathcal{C}^r$ maps  on a compact manifold with $r>1$ have symbolic extensions. The following theorem answers affirmatively to this conjecture in dimension $2$ and extends thus the results of \cite{burinv}.  When $T:M\rightarrow M$ is a $\mathcal{C}^1$ map on a compact Riemannian manifold $(M,\|\|)$ we denote by $R(T)$ the exponential growth of the derivative, i.e.  $R(T)=\lim_{n\rightarrow +\infty}\frac{\log^+\|DT^n\|}{n}$. This quantity does not depend on the choice of the Riemannian metric $\|\|$ on $M$.

\begin{theo}\label{intro}
Let $T:M\rightarrow M$ be a $\mathcal{C}^r$ map on a compact surface $M$ with $r>1$. Then $T$ admits symbolic extensions and

\begin{eqnarray*}h_{sex}(T)\leq h_{top}(T)+\frac{4R(T)}{r-1}
\end{eqnarray*}

Moreover, if $T$ is a local surface diffeomorphism, then

\begin{eqnarray*}h_{sex}(T)\leq h_{top}(T)+\frac{R(T)}{r-1}
\end{eqnarray*}

\end{theo}

The paper is organized as follows. We first recall the background of the theory of symbolic extensions and  properties of continuity of the sum of the positive Lyapunov exponents. Following S.Newhouse  we also recall how the local entropy is bounded from above by the local volume growth of smooth disks. Then we state our main results and as in \cite{burinv} we reduce them to a Reparametrization Lemma of Bowen's balls in a similar (but finer) approach of the classical Yomdin theory.  The last sections are devoted to the proof of the Reparametrization Lemma.

\section{Preliminaries}

In the following we denote $\mathcal{M}(X,T)$ the set of invariant Borel probability measures of the dynamical system $(X,T)$ and $\mathcal{M}_e(X,T)$ the subset of ergodic measures. We endow $\mathcal{M}(X,T)$ with the weak star topology. Since $X$ is a compact metric space, this topology is metrizable. We denote by $\dist$ a metric on $\mathcal{M}(X,T)$. It is well known that $\mathcal{M}(X,T)$ is  compact and  convex and its extreme points  are exactly the ergodic measures. Moreover if $\mu\in \mathcal{M}(X,T)$ there exists an unique Borel probability measure $M_{\mu}$ on $\mathcal{M} (X,T)$ supported by $\mathcal{M}_e(X,T)$ such that  for all Borel subsets $B$ of $X$ we have $\mu(B)=\int\nu(B)dM_{\mu}(\nu)$. This is the so called ergodic decomposition of $\mu$. A bounded real Borel map $f:\mathcal{M}(X,T)\rightarrow \R$ is said to be harmonic if $f(\mu)=\int_{\mathcal{M}_e(X,T)}f(\nu)dM_{\mu}(\nu)$ for all $\mu\in \mathcal{M}(X,T)$. It is a well known fact that affine upper semicontinuous maps are harmonic. The measure theoretical entropy $h:\mathcal{M}(X,T)\rightarrow \mathbb{R}^+$ is  always harmonic \cite{W} but is not upper semicontinuous in general. It may  not be upper semicontinuous even for $\mathcal{C}^r$ map for any $r\in \R^+$ \cite{Misss}. However $h$ is upper semicontinuous for $\mathcal{C}^{\infty}$ maps \cite{Neww}.

If $f$ is a bounded real Borel  map defined on $\mathcal{M}_e(X,T)$, the harmonic extension $\overline{f}$ of $f$ is the function defined on $\mathcal{M}(X,T)$ by :

$$\overline{f}(\mu):=\int_{\mathcal{M}_e(X,T)}f(\nu)dM_{\mu}(\nu)$$
It is easily seen that $\overline{f}$ coincides with $f$ on $\mathcal{M}_e(X,T)$ and that $\overline{f}$ is harmonic.

\subsection{Entropy structure}
The measure theoretical entropy function can be computed in many ways as limits of  nondecreasing sequences of  nonnegative functions defined on $\mathcal{M}(X,T)$ (with decreasing sequences of partitions, formula of Brin-Katok,...). The entropy structures are such particular
sequences whose convergence reflect the topological dynamic : they allow for example to compute the tail entropy \cite{Bur2} \cite{Dow2}, but also especially the symbolic  extension entropy \cite{BD} \cite{Dow2} (see below for precise statements).

We skip the formal definition of entropy structures, but we recall a basic fact. Two  nondecreasing sequences, $(h_k)_{k\in \N}$ and $(g_k)_{k\in\N}$, of nonnegative functions defined on $\mathcal{M}(X,T)$ are said to be uniformly equivalent  if for all $\gamma>0$ and for all  $k\in \N$, there exists $l\in\N$ such that $h_l>g_k-\gamma$ and $g_l>h_k-\gamma$. Two entropy structures are uniformly equivalent and any nondecreasing sequence of nonnegative functions which is uniformly equivalent to an entropy structure is itself an entropy structure. In other terms the set of entropy structures is an equivalence class for the above relation.

We recall now Lemma 1 of \cite{burinv} which relates  the entropy structures of a given dynamical system with those of its inverse and powers.

\begin{lem}\label{struun}
Let $(X,T)$ be a dynamical system with finite topological entropy and let $\mathcal{H}=(h_k)_{k\in \N}$ be an entropy structure of $T^p$ with $p\in \mathbb{N}\setminus \{0\}$ (when $T$ is a homeomorphism we consider $p\in \mathbb{Z}\setminus\{0\})$.  Then the sequence  $\frac{1}{|p|}\mathcal{H}|_{\mathcal{M}(X,T)}=\left(\frac{h_k|_{\mathcal{M}(X,T)}}{|p|}\right)_{k\in\N}$ defines an entropy structure of $T$.
\end{lem}

We finally check that the minimum of two entropy structures defines again an entropy structure.

\begin{lem}\label{strudeux}
Let $(X,T)$ be a dynamical system with finite topological entropy. If $\mathcal{H}=(h_k)_k$ and $\mathcal{G}=(g_k)_k$ are two entropy structures, then $\min(\mathcal{H},\mathcal{G}):=\left(\min(h_k,g_k)\right)_k$ is also an entropy structure.
\end{lem}

\begin{demo}
Let $\gamma>0$ and $k\in\N$. As $\mathcal{H}$ and $\mathcal{G}$ are both entropy structures, they are in particular uniformly equivalent. Therefore
there exists an integer $l$ such that  $h_l>g_k-\gamma$ and  $g_l>h_k-\gamma$. We can assume that $l>k$ so that $h_l\geq h_k$ by monotonicity of $\mathcal{H}$. Therefore $h_l>\min(h_k,g_k)-\gamma$ and  $\min(h_l,g_l)>h_k-\gamma$.
\end{demo}

\subsection{Tail entropy}
In order to study the properties of upper semicontinuity of the entropy function of a dynamical system and in particular the existence of measures of maximal entropy, M.Misiurewicz introduced in the seventies the following quantity which is now known as the tail entropy of the system. Let us first recall some usual notions relating to the entropy of dynamical systems (we refer to \cite{W} for a general introduction to entropy).

 Consider a continuous map $T:X\rightarrow X$ with $(X,d)$  a compact metric space. Let $n\in \N$ and $\delta>0$. A subset $F$ of $X$ is called a $(n,\delta)$ separated set when for all $x,y\in F$ there exists $0\leq k< n$ such that $d(T^kx,T^ky)\geq\delta$. Let $Y$ be a subset of $X$. A subset $F$ of $Y$ is called a $(n,\delta)$ spanning  set of $Y$ when for all $y\in Y$ there exists $z\in F$ such that $d(T^ky,T^kz)<\delta$ for all $0\leq k< n$.
 Given a point $x\in X$ we denote by $B(x,n,\delta)$ the Bowen's ball centered at $x$ of radius $\delta$ and length $n$ :
$$B(x,n,\delta):=\{y\in X, \ d(T^ky,T^kx)<\delta \ \text{for} \ k=0,...,n-1\}$$

The tail entropy, $h^*(T)$,  of $(X,T)$ is then defined by
$$h^{*}(T):=\lim_{\epsilon\rightarrow 0}\limsup_{n\rightarrow +\infty}\sup_{x\in X}\frac{1}{n}\log \min \left\{\sharp C \ : \ F \text{ is a } (n,\delta) \text{ spanning set of } B(x,n,\epsilon)\right\} $$

This quantity is a topological invariant which estimates the entropy appearing at arbitrarily small scales.
The tail entropy bounds from above the defect of upper semicontinuity of the entropy function \cite{Mis} :

\begin{eqnarray*}
\forall \mu\in \mathcal{M}(X,T), \ \limsup_{\nu\rightarrow \mu}h(\nu)-h(\mu)\leq h^*(T)
\end{eqnarray*}

 In general the supremum of the defect of upper semicontinuity of the entropy function differs from the tail entropy. But it is easily seen that for any entropy structure $(h_k)_k$ of $(X,T)$, we have $\limsup_{\nu\rightarrow \mu}h(\nu)-h(\mu)\leq \lim_k\limsup_{\nu\rightarrow \mu}(h-h_k)(\nu)$ and T.Downarowicz proved then the following variational principle \cite{Dow2} (see also \cite{Bur2}) :

\begin{eqnarray}\label{tailvar}
\sup_{\mu\in \mathcal{M}(X,T)}\lim_k\limsup_{\nu\rightarrow \mu}(h-h_k)(\nu)=\lim_k\sup_{\mu\in \mathcal{M}(X,T)}(h-h_k)(\mu)=h^*(T)
\end{eqnarray}

By using Yomdin theory J.Buzzi \cite{Buz} established the following upper bound on the tail entropy of $\mathcal{C}^r$ maps $T$ on a compact manifold $M$ with $r>1$ :

\begin{eqnarray}\label{buz}
h^*(T)\leq \frac{\dim(M)}{r}R(T)\end{eqnarray}

This inequality is known to be sharp for noninvertible maps \cite{Buz}, \cite{Rue}. We will prove in the present paper a similar sharp upper bound on the tail entropy of $\mathcal{C}^r$ surface diffeomorphisms with $r>1$ (see Theorem \ref{tail} below).

When $h^*(T)=0$ the dynamical system $(X,T)$  is said to be asymptotically $h$-expansive. For example, uniformly hyperbolic dynamical systems or piecewise monotone interval maps are asymptotically $h$-expansive. Then entropy structures are converging uniformly to the entropy function. There exist therefore measures of maximal entropy by upper semicontinuity of the entropy function and  also symbolic extensions preserving the entropy of invariant measures (see Theorem \ref{nouveau} below).  According to Inequality (\ref{buz}) it is also the case of $\mathcal{C}^{\infty}$ maps on a compact manifold.

\subsection{Symbolic extension entropy function}
 A symbolic extension of $(X,T)$ is a subshift $(Y,S)$ of a full shift on a finite number of symbols, along with a continuous surjection $\pi: Y \rightarrow X$ such that $T\circ  \pi = \pi\circ S$. Given a symbolic extension $\pi:(Y,S)\rightarrow(X,T)$ we consider the extension entropy  $h^{\pi}_{ext}:\mathcal{M}(X,T)\rightarrow \R^+$ defined for all $\mu\in \mathcal{M}(X,T)$ by :
$$h^{\pi}_{ext}(\mu)=\sup_{\pi^*\nu=\mu}h(\nu)$$

Then the symbolic extension entropy function, $h_{sex}:\mathcal{M}(X,T)\rightarrow \R^+$,  is :
$$h_{sex}=\inf h_{ext}^{\pi}$$

where the infimum holds over all the symbolic extensions of $(X,T)$. By convention, if $(X,T)$ does not admit any symbolic extension we simply put $h_{sex}\equiv+\infty$. Recall we defined in the introduction the topological symbolic extension entropy $h_{sex}(T)$ as the infimum of the topological entropy of the symbolic extensions of $(X,T)$ (as previously we put $h_{sex}(T)=+\infty$ when there are no such extensions). M.Boyle and T.Downarowicz proved that these two quantities are related by the following variational principle :

\begin{eqnarray}\label{pvsex}
h_{sex}(T)&=&\sup_{\mu\in \mathcal{M}(X,T)}h_{sex}(\mu)
\end{eqnarray}

We present now the major Symbolic Extension Entropy Theorem of \cite{BD} which allows to compute the symbolic extension entropy function from the properties of convergence of any entropy structure. We follow the exposition in \cite{buren}.
 Let $\mathcal{S}(X,T)$ be the set of nonnegative upper semicontinuous functions defined on $\mathcal{M}(X,T)$ to which we add the function
 constant equal to $+\infty$. Let $(h_k)_k$ be an entropy structure, we define an increasing  operator  on $\mathcal{S}(X,T)$, denoted by $\mathcal{T}_{sex}$,  as follows :

\begin{eqnarray*}
\mathcal{T}_{sex}: \mathcal{S}(X,T) &\rightarrow & \mathcal{S}(X,T)\\
f&\mapsto & \left[\mu\mapsto \lim_{k} \limsup_{\nu\in \mathcal{M}(X,T), \nu \rightarrow \mu}\left(f+h-h_k\right)(\nu)\right]
\end{eqnarray*}

One easily checks from the uniform equivalence relation that $\mathcal{T}_{sex}$ does not depend on the choice of the entropy structure $(h_k)_k$. We also observe that the tail variational principle can be written as $\sup_{\mu\in \mathcal{M}(X,T)}\mathcal{T}_{sex}0(\mu)=h^*(T)$.
By using the affine structure of the set of invariant probability measures M.Boyle and T.Downarowicz proved that the least fixed point of $\mathcal{T}_{sex}$ (whose existence is ensured by Tarski-Knaster Theorem) coincides with the infimum of the affine fixed points of $\mathcal{T}_{sex}$. We recall now the Symbolic Extension Entropy Theorem :

\begin{theo}(Theorem 5.5 of  \cite{BD})\label{nouveau}
The affine fixed points of $\mathcal{T}_{sex}$ are exactly the functions $h_{ext}^{\pi}-h$, i.e. $f$ is a nonnegative affine upper semicontinuous function
on $\mathcal{M}(X,T)$ fixed by $\mathcal{T}_{sex}$ if and only if  there exists a symbolic extension $\pi$ such that $f=h^{\pi}_{ext}-h$. Moreover $h_{sex}-h$ is the least fixed point of $\mathcal{T}_{sex}$.
\end{theo}

When $(X,T)$ is asymptotically $h$-expansive, then the zero function is a fixed point of $\mathcal{T}_{sex}$ because  entropy structures are converging uniformly. Therefore such systems (including $\mathcal{C}^{\infty}$ maps on compact manifold) admit symbolic extensions $\pi:(Y,S)\rightarrow (X,T)$ with $h^{\pi}_{ext}=h$. This result was first proved by M.Boyle, D.Fiebig and U.Fiebig \cite{bff}.\\

The Passage Theorem of T.Downarowicz and A.Maass \cite{Dow} gives a weaker condition for an affine upper semicontinuous function
on $\mathcal{M}(X,T)$ to be fixed by $\mathcal{T}_{sex}$ : one only needs to consider the $\limsup$ in ergodic measures $\nu$ in the definition of $\mathcal{T}_{sex}$. More precisely a nonnegative affine upper semicontinuous function $f$
on $\mathcal{M}(X,T)$   is a fixed point of $\mathcal{T}_{sex}$ if and only if  $\limsup_{\nu\in \mathcal{M}_e(X,T), \nu \rightarrow \mu}\left(f+h-h_k\right)(\nu)=f(\mu)$ for all $\mu\in \mathcal{M}(X,T)$. This can be restated as follows :

\begin{theo}[Downarowicz, Maass]\cite{Dow}\cite{dowbook}\label{theopass}
Let $(X,T)$ be a dynamical system with finite topological entropy.  Let $g$ be a nonnegative upper semicontinuous affine function on $\mathcal{M}(X,T)$ such that for every $\gamma>0$, for every $\mu\in \mathcal{M}(X,T)$ and for every entropy structure $\mathcal{H}=(h_k)_{k\in\N}$ there exists  $k_{\mu}\in\N$ and $\delta_{\mu}>0$ such that for every ergodic measure $\nu$ satisfying $\dist(\nu,\mu)<\delta_{\mu}$ it holds that :

\begin{eqnarray}\label{estim}(h-h_{k_{\mu}})(\nu)\leq g(\mu)-g(\nu)+\gamma\end{eqnarray}
Then there exists a symbolic extension $\pi:(Y,S)\rightarrow (X,T)$ such that $h^{\pi}_{ext}-h=g$. In particular $h_{sex}-h\leq g$.

\end{theo}

It follows from the uniform equivalence relation that  the assumptions of Theorem \ref{theopass} hold as
soon as they are satisfied by one particular entropy structure.

\subsection{Newhouse local entropy}\label{newent}

We recall now the "Newhouse local entropy". Let $x\in X$, $\epsilon>0$, $\delta>0$, $n\in \N$ and $F\subset X$ a Borel set, we define :

$$H(n,\delta|x,F,\epsilon):=\log \max\left\{\sharp E \ : \  E\subset F\bigcap B(x,n,\epsilon) \ \text{and} \  E \ \text{is a}  \ (n,\delta) \ \text{separated set}\right\} $$

$$H(n,\delta|F,\epsilon):=\sup_{x\in F} H(n,\delta|x,F,\epsilon)$$

$$h(\delta|F,\epsilon):=\limsup_{n\rightarrow +\infty}\frac{1}{n}H(n,\delta|F,\epsilon)$$

$$h(X|F,\epsilon):=\lim_{\delta\rightarrow 0}h(\delta|F,\epsilon)$$

 Then for any ergodic measure $\nu$ we put :
$$h^{New}(X|\nu,\epsilon):=\lim_{\alpha\rightarrow 1}\inf_{\nu(F)>\alpha}h(X|F,\epsilon)$$

Finally we extend the function $h^{New}(X|\cdot,\epsilon)$ to $\mathcal{M}(X,T)$ by the harmonic extension. Given a nonincreasing sequence $(\epsilon_k)_{k\in \N}$ converging to $0$, we consider the sequence
$\mathcal{H}_T^{New}=(h_k^{New})_{k\in\N}$ with $h_k^{New}:=h-h^{New}(X|.,\epsilon_k)$ for all integers $k$. A similar quantity with minor differences  was first introduced by S.Newhouse in \cite{Neww}. T.Downarowicz proved this sequence defines an entropy structure \cite{Dow2} for homeomorphisms and the author extended the result in the noninvertible case \cite{Bur}.

\subsection{Lyapunov exponents}

Let $(M,\|\|)$ be a compact Riemannian manifold of dimension $d$ and let $T:M\rightarrow M$ be a $\mathcal{C}^1$ map.  We denote $\|D_xT\|=\sup_{v\in T_xM-\{0\}}\frac{\|D_xT(v)\|}{\|v\|}$ the induced norm of the differential $D_xT$ of $T$ at $x$  and  $\|DT\|=\sup_{x\in M}\|D_xT\|$ the supremum norm of the differential of $T$. We consider an ergodic $T$-invariant measure $\nu$. According to Oseledet's theorem \cite{Osc}, there exist real numbers  $\chi_1(\nu)\geq ...
\geq \chi_d(\nu)\geq -\infty$, an increasing  sequence of measurable invariant subbundles of the tangent space  $\{0\}=G_{d+1}\subseteq G_d \subseteq...\subseteq G_1=TM$  and a Borel set $F$
with $\nu(F)=1$ such that for all $x\in F$ and all $v_i\in G_{i}\setminus G_{i+1}$ with $i=1,...,d$ we have
\begin{eqnarray*}
  \lim_{n\rightarrow +\infty}\frac{1}{n}\log\|D_xT^n(v_i)\|&=&\chi_{i}(\nu)
\end{eqnarray*}

 The real numbers $\chi_i(\nu)$ are the well-known Lyapunov exponents of $\nu$ (sometimes we use also the notations $\chi_i(\nu,T)$  to avoid ambiguities). We prove now elementarily that the harmonic extension of the sum  of the positive Lyapunov exponent is upper semicontinuous.  In the one dimensional case it was proved by T.Downarowicz and A.Maass by using a clever argument of convexity (see Fact 2.5 of \cite{Dow}). We just adapt the proof of \cite{burinv} which only deals with the dimension $d=2$. If $S:M\rightarrow M$ is a $\mathcal{C}^1$ map we denote by $\Lambda^kD_xS$ the map induced by the differential map $D_xS$ on the $k^{th}$ exterior power $\Lambda^kT_xM$ (with $1\leq k\leq d$) of the tangent space of $M$ at $x\in M$ and by $\|\|_k$ the induced norm on the space of multilinear maps on $\Lambda^kT_xM$.

\begin{lem}\label{thes}
For all $1\leq e\leq d$ and for all $\mu\in \mathcal{M}(M,T)$, we have :

$$\overline{\sum_{i=1}^e\chi_i^+}(\mu)=\inf_{n\in\N}\frac{1}{n}\int \max_{k=1,...,e}\log^+\|\Lambda^kD_xT^n\|_kd\mu(x)$$

 In particular $\overline{\sum_{i=1}^e\chi_i^+}:\mathcal{M}(M,T)\rightarrow \R^+$ is upper semicontinuous.
\end{lem}

\begin{demo}

For all integers $n>0$ we consider the function $f_n:\mathcal{M}(M,T)\rightarrow \R^+$ defined by :
$$\forall \mu \in \mathcal{M}(M,T), \ f_n(\mu)=\int \max_{k=1,...,e}\log^+\|\Lambda^kD_xT^n\|_kd\mu(x)$$
This function is clearly continuous and affine, and therefore harmonic. Also  $(f_n(\mu))_{n\in\N}$ is a subadditive sequence for all $\mu\in \mathcal{M}(M,T)$.

It follows from Oseledet's theorem that $\sum_{i=1}^e\chi_i^+(\nu)=\lim_{n\rightarrow +\infty}\frac{f_n(\nu)}{n}$ for all ergodic measures $\nu$.
Consider now a general measure $\mu\in \mathcal{M}(M,T)$. We have :

\begin{eqnarray*}
\overline{\sum_{i=1}^e\chi_i^+}(\mu)&:=&\int_{\mathcal{M}_e(M,T)}\sum_{i=1}^e\chi_i^+(\nu)dM_{\mu}(\nu)\\
 &=&\int_{\mathcal{M}_e(M,T)}\lim_{n\rightarrow +\infty}\frac{f_n(\nu)}{n}dM_{\mu}(\nu)
\end{eqnarray*}

Obviously $f_n(\nu)\leq e\log^+\|DT\|$ for all ergodic measures $\nu$. Therefore by applying the theorem of dominated convergence we get :

\begin{eqnarray*}
\overline{\sum_{i=1}^e\chi_i^+}(\mu)&=&\lim_{n\rightarrow +\infty}\int_{\mathcal{M}_e(M,T)}\frac{f_n(\nu)}{n}dM_{\mu}(\nu)
\end{eqnarray*}

and by harmonicity of $f_n$ :

\begin{eqnarray*}
\overline{\sum_{i=1}^e\chi_i^+}(\mu)&=&\lim_{n\rightarrow +\infty}\frac{f_n(\mu)}{n}
\end{eqnarray*}

But the sequence $(f_n(\mu))_{n\in \N}$ is subadditive so that :

\begin{eqnarray*}
\overline{\sum_{i=1}^e\chi_i^+}(\mu)&=&\inf_{n\in\N}\frac{f_n(\mu)}{n}
\end{eqnarray*}

We conclude that $\overline{\sum_{i=1}^e\chi_i^+}$ is an upper semicontinuous function as an infimum of a family of continuous functions.
\end{demo}

Clearly the function $\overline{\sum_{i=1}^e\chi_i^+}$ is uniformly bounded from above by
$$R_e(T):=\limsup_n \sup_{x\in X}\frac{1}{n}\max_{k=1,...,e}\log^+\|\Lambda^kD_xT^n\|_k$$ With the previous notations we have also
$R_e(T)\leq e R(T)$ and $R_1(T)=R(T)$. In fact  the following variational principle holds \cite{Bur} :

$$\sup_{\mu}\overline{\sum_{i=1}^e\chi_i^+}(\mu)=R_e(T)$$

We notice that the harmonic extension $\overline{\sum_{i=1}^e\chi_i^+}$ can be rewritten in the most common way as
$\overline{\sum_{i=1}^e\chi_i^+}=\int_M\sum_{i=1}^e \chi_i^+(x)d\mu(x)$ where  $(\chi_i(x))_{i=1,...,d}$ denotes the Lyapunov exponents of a regular point $x\in M$. \\

In the special case $e=d$ we also recall that the sum of all the  Lyapunov exponents (positive or not) are given by

$$\overline{\sum_{i=1}^d\chi_i}(\mu)=\int_M \log\Jac_x(T)d\mu(x)$$

In the following we are interesting in the entropy of ergodic measures.  We denote by $\sum\chi^+$ (resp. $\sum\chi^-$) the sum of all the positive  (resp. negative)  Lyapunov exponents. The  Ruelle-Margulis inequality states that for a $\mathcal{C}^1$ map $T:M\rightarrow M$ on a compact manifold $M$ the entropy $h_T(\nu)$ of an ergodic measure $\nu$ is bounded from above by the sum of its positive Lyapunov exponents. When $T$ is a  diffeomorphism it is easily seen  by applying the Ruelle-Margulis inequality to $T$
 and its inverse $T^{-1}$ and by using the equality $h_T(\nu)=h_{T^{-1}}(\nu)$ that  any ergodic measure $\nu\in \mathcal{M}_e(M,T)$ with non zero entropy has  at least one positive and one negative Lyapunov exponent,  moreover

$$h(\nu)\leq \min\left(\sum\chi^+(\nu),-\sum\chi^-(\nu)\right)$$

\subsection{Local entropy and volume growth}\label{volume}
In this section we recall Theorem 2 of \cite{Neww} (in fact an intermediate result in its proof) which relates the Newhouse local entropy of an ergodic measure with the local volume growth of smooth disks. We begin with some definitions. Let $T:M\rightarrow M$ be a $\mathcal{C}^r$ map on a compact manifold with $r>1$. We fix a Riemannian metric $\|\|$ on the manifold $M$ and we endow $M$ with the induced distance.

A $\mathcal{C}^1$ map $\sigma$ from the unit square  $[0,1]^k$ of $\R^k$ to $M$, which is a diffeomorphism onto its image, is called a  $k$-disk. For any $k$-disk $\sigma$ and for any $\chi>0$, $\gamma>0$ and $C>1$, we consider the set $\mathcal{H}^n_T(\sigma,\chi,\gamma,C)$ of points of $[0,1]^k$ whose exponential growth of the derivative of the $n$-first iterations of $T$ composed with $\sigma$ is almost equal to $\chi$ :

$$\mathcal{H}^n_T(\sigma,\chi,\gamma,C):=\left\{t\in [0,1]^k \ :  \ \forall 1\leq j\leq n, \ C^{-1}e^{(\chi-\gamma)j}\leq \| D_t\left(T^j\circ\sigma\right)\|\leq Ce^{(\chi+\gamma)j}\right\}$$

 We also denote by $|\sigma|$ the $k$-volume of $\sigma$, i.e. $|\sigma|=\int|\Lambda^kD_t\sigma|d\lambda(t)$ where $d\lambda$ is the Lebesgue measure on $[0,1]^k$. Then given  $\chi>0$, $\gamma>0$, $C>1$, $x\in M$, $n\in \N$ and $\epsilon >0$ we define the local  volume growth of $\sigma$ at $x$ with respect to these parameters as follows :

\begin{eqnarray*} V_x^{n,\epsilon}(\sigma|\chi,\gamma,C):= \left| T^{n-1}\circ \sigma|_{\mathcal{H}^n_T(\sigma,\chi,\gamma,C) \cap \sigma^{-1}\left(B(x,n,\epsilon)\right)}\right|
\end{eqnarray*}

Fix an ergodic measure $\nu$ and let $l_{\nu}$ be  the codimension  of the first Lyapunov eigenspace $G_i$ with $\chi_i(\nu)<0$. In \cite{Neww} S.Newhouse proved by using Pesin theory that for some $\epsilon>0$ depending only on the manifold $M$ and for any $1>\alpha>0$ and  $\gamma>0$  there exist a Borel set $F_\alpha$ of $\nu$-measure larger than $\alpha$ and a $l_{\nu}$-disk $\sigma$ of class $\mathcal{C}^{\infty}$, such that for all $\delta>0$ one can find a constant $C>1$ depending only on $F_\alpha$ and another constant $D$ depending only on  $\delta$ and $F_{\alpha}$ satisfying for all $n\in \N$ and for all  $x\in F_{\alpha}$  :

\begin{eqnarray}\label{vvolume}H(n,\delta|x,F_{\alpha},\epsilon)\leq De^{\gamma n}V_x^{n,\epsilon}(\sigma|\chi_1(\nu),\gamma,C)
\end{eqnarray}

  In fact $F_{\alpha}$ is a Pesin set and  the $l_{\nu}$-disk can be chosen to be the exponential map restricted to a neighborhood of the unstable distribution at some point of $F_{\alpha}$. In the following we will consider only smooth $1$-disks $\sigma:[0,1]\rightarrow \R^d$. We  denote by $\sigma'$ the derivative $D\sigma$ of the curve $\sigma$.

\section{Statements}
We first state in this section the main result of this paper by specifying Theorem \ref{intro}  at the measure theoretic level. It will be deduced, by using the Estimate theorem, from the below Main Proposition which gives an upperbound on the Newhouse local entropy involving  the positive Lyapunov exponents. The proof of this proposition is the topic of the next sections.

\begin{theo}\label{corps}
Let  $T:M\rightarrow M$ be a $\mathcal{C}^r$  surface map (resp. local surface diffeomorphism) with $r>1$. Then there
exists a symbolic extension $\pi:(Y,S)\rightarrow (X,T)$ with

$$h^{\pi}_{ext}-h=\frac{2\overline{\sum\chi^+}(\mu)}{r-1} \ \left(\text{resp. } =\frac{\overline{\chi_1^+}(\mu)}{r-1}\right)$$

 In particular, for all $\mu\in \mathcal{M}(M,T)$,

$$h_{sex}(\mu)- h(\mu)\leq\frac{2\overline{\sum\chi^+}(\mu)}{r-1} \  \left(\text{resp. }\frac{\overline{\chi_1^+}(\mu)}{r-1}\right)$$
\end{theo}

The above theorem  will follow  from the following Main Proposition by applying the Estimate Theorem to the upper semicontinuous affine function $g=\frac{2\overline{\sum\chi^+}}{r-1}$ or $g=\frac{\overline{\chi_1^+}}{r-1}$. Remark  that Theorem \ref{intro} stated in the introduction is the topological version of  Theorem \ref{corps} : it is deduced by the usual variational principle for the entropy and the variational principle for the symbolic extension entropy (Equation (\ref{pvsex})). According to the examples built in \cite{ND}  the upper bound on the symbolic extension entropy function in Theorem \ref{corps} is sharp in the case of local surface diffeomorphism. More precisely S.Newhouse and T.Downarowicz exhibit examples of $\mathcal{C}^r$ surface diffeomorphisms with $r>1$ which admit an hyperbolic periodic measure $\gamma_p$ with $h_{sex}(\gamma_p)=\frac{\chi_1^+(\gamma_p)}{r-1}$. However the optimality of the upper bound in the general case is open.

As the functions $h^{\pi}_{ext}$, $\overline{\chi_1^+}$ and $\overline{\sum\chi^+}$ are upper semicontinuous it follows from Theorem \ref{corps} that the entropy function of a $\mathcal{C}^r$ surface map with $r>1$ is a difference of nonnegative upper semicontinuous functions. In particular generic measures are continuity points of the entropy function. In general this is false for $\mathcal{C}^1$ maps. We refer to \cite{buren} for further details on the links between the symbolic extension entropy and the properties of continuity of the entropy function.\\

We give now a new upper bound on the tail entropy of $\mathcal{C}^r$ surface diffeomorphisms with $r>1$. The examples built in \cite{ND}
show that this upper bound is sharp.

\begin{theo}\label{tail}
Let $T:M\rightarrow M$ be a $\mathcal{C}^r$ surface diffeomorphisms with $r>1$. Then, for all $\mu\in \mathcal{M}(X,T)$, we have
$\mathcal{T}_{sex}0(\mu)\leq \frac{\overline{\chi_1^+}(\mu)}{r}$ and in particular,

$$h^*(T)\leq \frac{R(T)}{r}$$
\end{theo}

 Theorem \ref{corps} and Theorem \ref{tail} will both follow  from the following proposition.

\begin{Prop}\label{gg}
Let $T:M\rightarrow M$ be a $\mathcal{C}^r$ map with $r>1$ on a compact manifold $M$ of dimension $d$. Let $\mu$ be a $T$-invariant measure and fix some $\gamma> 0$.\\

 Then there exist $\delta_{\mu}>0$, an entropy structure $(h_k)_k$ and $k_{\mu}\in \N$ such that for
every ergodic $T$-invariant measure $\nu$ with $\dist(\nu,\mu) < \delta_{\mu}$ and  either with at most one nonnegative Lyapunov exponent or with all Lyapunov exponents nonnegative (i.e.   $l_{\nu}$  is equal to either $0$, $1$ or $d$) it holds that :

\begin{eqnarray}\label{ggeq}
(h-h_{k_{\mu}})(\nu)\leq \frac{l_{\nu}}{r-1}\left(\overline{\sum_{i=1}^{l_{\nu}}\chi^+}(\mu)-\sum\chi^+(\nu)\right)+\gamma
\end{eqnarray}

\end{Prop}

In the case $l_{\nu}=d$ it is a direct consequence of the Main Theorem of \cite{superbur} that we recall below. Then, to prove the proposition we only need to consider the case $l_{\nu}\leq 1$ and we finally conclude by considering the minimum of the two entropy structures, the one corresponding to the case $l_{\nu}\leq1$ and the other corresponding to the case $l_{\nu}=d$ (see Lemma \ref{strudeux}). It is worth noting that when $d=1$ Inequality (\ref{ggeq}) can be obtained by two ways which seem to be different : the first one due to T.Downarowicz and A.Maass by studying the critical set and the second one presented in the present paper which follows from the Reparametrization Lemma stated in the next section.

We think that this statement should hold for any ergodic measure without condition on its Lyapunov spectrum. It will implies Theorem \ref{corps} in any dimension. Let us just do it in the two dimensional case. We first recall the main result of \cite{superbur} which is a  generalization of the approach of T.Downarowicz and A.Maass.

\begin{prop}(Main Theorem of \cite{superbur})
Let $T:M\rightarrow M$ be a $\mathcal{C}^r$ map, with $r>1$,  on a compact manifold of dimension $d$.  Let $\mu$ be an invariant measure and fix some $\gamma> 0$.\\

Then there exist $\delta_{\mu}>0$, an entropy structure $(h_k)_k$ and $k_{\mu}\in \N$  such
 that for every ergodic measure $\nu$ with $\dist(\nu,\mu)< \delta_{\mu}$ it holds that :

\begin{equation}\label{tuer}
(h-h_{k_\mu})(\nu)\leq \frac{d\left(f(\mu)-f(\nu)\right)}{r-1}-\sum\chi^-(\nu)+\gamma
\end{equation}
where $f(\xi)=\max\left(\int\log \Jac _x(T)d\xi(x),0\right)$ for all invariant measures $\xi$.
\end{prop}

When $T$ is a local diffeomorphism then $\xi\mapsto \int\log \Jac _x(T)d\xi(x)$ is continuous on $\mathcal{M}(M,T)$. Then, with the notations of the above theorem, we have  $(h-h_{k_\mu})(\nu)\leq -\sum\chi^-(\nu)+2\gamma$ for $\nu$ close enough to $\mu$. In particular if $(h-h_{k_\mu})(\nu)>2\gamma$ then $\nu$ has at least one negative Lyapunov exponent (and also at least a positive one by Ruelle-Margulis inequality).\\

We deduce now the theorems of this section from the above proposition.

\begin{demof}{Theorem \ref{corps}}
As we follows the strategy of \cite{Dow} we only sketch the proof. First we consider the case of local surface diffeomorphisms. According to the above remark it is enough to consider ergodic measures  with one negative and one positive Lyapunov exponent. Then one just applies the Estimate Theorem with $g=\frac{\overline{\chi_1^+}}{r-1}$ which is an upper semicontinuous function satisfying (\ref{estim}) by assumption. We prove similarly Theorem \ref{corps} for general surface maps by applying the Estimate Theorem with $g=\frac{2\overline{\sum\chi^+}}{r-1}$.

\end{demof}

\begin{demof}{Theorem \ref{tail}}
Let $(h_k)_k$ be an entropy structure of $(M,T)$ and let $\mu \in \mathcal{M}(X,T)$. By Ruelle-Margulis inequality we have for all $\nu\in \mathcal{M}(M,T)$

$$(h-h_k)(\nu)\leq \min\left(\overline{\chi_1^+}(\nu),\left(h-h_k+\frac{\overline{\chi_1^+}}{r-1}\right)(\nu)-\frac{\overline{\chi_1^+}(\nu)}{r-1}\right)$$

By taking the $\limsup$ when $\nu$  goes to $\mu$ we have for all integers $k$

\begin{eqnarray*}
\limsup_{\nu\rightarrow \mu}(h-h_k)(\nu)&\leq & \sup_{\nu\in \mathcal{M}(X,T)}\min\left(\overline{\chi_1^+}(\nu), \limsup_{\xi\rightarrow \mu}\left(h-h_k+\frac{\overline{\chi_1^+}}{r-1}\right)(\xi)-\frac{\overline{\chi_1^+}(\nu)}{r-1}\right)
\end{eqnarray*}

and  since $\frac{\overline{\chi_1^+}}{r-1}$ is a fixed point of $\mathcal{T}_{sex}$ by Theorem \ref{corps}  we have when $k$ goes to infinity

\begin{eqnarray*}
\lim_k\limsup_{\nu\rightarrow \mu}(h-h_k)(\nu)&\leq & \sup_{\nu\in \mathcal{M}(X,T)}\min\left(\overline{\chi_1^+}(\nu), \lim_k\limsup_{\xi \rightarrow \mu}\left(h-h_k+\frac{\overline{\chi_1^+}}{r-1}\right)(\xi)-\frac{\overline{\chi_1^+}(\nu)}{r-1}\right)\\
&\leq & \sup_{\nu\in \mathcal{M}(X,T)}\min\left(\overline{\chi_1^+}(\nu), \frac{\overline{\chi_1^+}(\mu)}{r-1}-\frac{\overline{\chi_1^+}(\nu)}{r-1}\right)
\end{eqnarray*}

As a continuous piecewise affine function in $\overline{\chi_1^+}(\nu)$ the right member $\min\left(\overline{\chi_1^+}(\nu), \frac{\overline{\chi_1^+}(\mu)}{r-1}-\frac{\overline{\chi_1^+}(\nu)}{r-1}\right)$ attains its maximum at $\frac{\overline{\chi_1^+}(\mu)}{r}$ and is therefore bounded from above by $\frac{\overline{\chi_1^+}(\mu)}{r}$. We conclude according to the tail variational principle (\ref{tailvar}) that $h^*(T)\leq \frac{R(T)}{r}$.
\end{demof}

We already have noticed that the inequality $h^*(T)\leq \frac{dR(T)}{r}$ due to J.Buzzi holds for any $\mathcal{C}^r$ maps on a compact manifold $M$ of dimension $d$. By adapting the proof it is not difficult to also prove $\mathcal{T}_{sex}0(\mu)\leq \frac{d\overline{\chi_1^+}(\mu)}{r}$ for all
$\mu\in \mathcal{M}(M,T)$ \cite{Bur}. But we do not know if in this inequality we can replace $dR(T)$ (resp. $d\overline{\chi_1^+}(\mu)$) by $\max_{1\leq e\leq d}R_e(T)$ (resp. $\overline{\sum\chi^+}(\mu)$).

\section{Bounding local entropy of curves with a Reparametrization Lemma of Bowen balls}

  In a similar way as in \cite{burinv} we reduce now the Main Proposition to a result of reparametrization of $1$-disk $\sigma$ in Bowen balls  by volume contracting maps. Then the local volume growth and thus the Newhouse local entropy of some ergodic measure $\nu$ is just bounded from above by the logarithmic growth of the number of reparametrizations. In fact we only need to reparametrize Bowen balls on a set of large $\nu$ measure. As one can assume the entropy and thus the maximal Lyapunov exponent $\chi_1(\nu)$ to be nonzero it is enough to consider the intersection of Bowen balls with the set $\mathcal{H}^n_T(\sigma,\chi_1(\nu),\gamma,C)$ for small $\gamma>0$ and for any $C>1$. This is one of the main differences with Yomdin theory where the whole Bowen ball is reparametrized by  volume contracting maps.  \\

In \cite{burinv} we do not use Newhouse estimate (which involves Pesin theory) of the local entropy with the local volume growth. Then we had to reparametrize not only curves but the  intersection of Bowen balls with \emph{finite hyperbolic sets} which are in general not one dimensional. The  present situation is of course easier and allows us to get $\mathcal{C}^r$ estimates.\\

 We denote $H:[1,+\infty[\rightarrow \R$ the function defined by $H(t)=-\frac{1}{t}\log(\frac{1}{t})-(1-\frac{1}{t})\log(1-\frac{1}{t})$. Moreover $[x]$ is the usual integer part of $x$ if  $x>0$ and zero if not (i.e. $[x]$ is the largest nonnegative integer $k$ such that $\max(x,0)\geq k$). Afterwards we also use the notation $\left\lceil x\right\rceil$ to denote the usual ceiling function (i.e. $\left\lceil x\right\rceil$ is the smallest integer $k$ such that $x\leq k$).  We recall that a map $T$ between two smooth Riemannian manifolds is said  of class $\mathcal{C}^r$ with $r>0$ when it admits a derivative of order $\left\lceil r-1\right\rceil$ which is $r-\left\lceil r-1\right\rceil$-Hölder. We denote then by $\|T\|_r$ the $r-\left\lceil r-1\right\rceil$-Hölder norm of the $\left\lceil r-1\right\rceil$-derivative of $T$. Observe that if $T$ is $\mathcal{C}^r$ then $T$ is $\mathcal{C}^s$ for all
 $0\leq s\leq r$. When $\left\lceil r-1\right\rceil < s\leq r$ we have moreover $\|T\|_s\leq C\|T\|_r$ with a constant $C$ depending only on the Riemannian manifolds. Finally $\|T\|_0$ will denote the usual supremum norm of $T$.

\begin{Lem}\label{main}
Let $T:M\rightarrow M$ be a $\mathcal{C}^r$ map with $r>1$ on a compact Riemannian manifold $M$.\\

  Then for all $\chi>0$, $\gamma>0$ and $C>1$,   there exist $\epsilon>0$  depending only on  $\|T\|_s$, $s=1,...,[r],r$ and a universal constant $A>0$  with the following properties. \\

 For all $x\in M$, for all positive integers $n$ and for all  $1$-disks  $\sigma:[0,1]\rightarrow M$ of class $\mathcal{C}^r$ with $\max_{s=1,...,[r],r}\|\sigma\|_s\leq 1$, there exist a family $\mathcal{F}_n$ of affine maps from $[0,1]$ to $[0,1]$ and a real positive number $B$ depending only on  $\sigma$, such that  with $\lambda^+_n(x,T):=\frac{1}{n}\sum_{j=0}^{n-1}\log^+\|D_{T^jx}T\|$ the following properties hold :

\begin{enumerate}[(i)]
\item $$\forall \psi\in \mathcal{F}_n, \ \forall 0 \leq l\leq n, \ \|D(T^{l}\circ \sigma\circ \psi)\|\leq 1,$$
\item $$\mathcal{H}^n_T(\sigma,\chi,\gamma,C)\cap \sigma^{-1}\left(B(x,n+1,\epsilon)\right)\subset \bigcup_{\psi\in \mathcal{F}_n}\psi([0,1]),$$
\item $$\log\sharp \mathcal{F}_n \leq  \frac{1}{r-1}\left(1+H([\lambda^+_n(x,T)-\chi]+3)\right)\left(\lambda^+_n(x,T)-\chi\right)n+An+B.$$
\end{enumerate}
\end{Lem}


We  deduce now the Main  Proposition  from the above statement by following \cite{burinv}. In the proof the terms $\lambda^+_n(x,T)$ for typical $\nu$ points $x$ and $\chi$ will be respectively related with the maximal Lyapunov exponent of $\mu$ and $\nu$ where $\nu$ is an ergodic measure near an invariant measure $\mu$. Moreover the quantity  $H([\lambda^+_n(x,T)-\chi]+3)$ will be negligible.

\begin{demof}{the Main Proposition assuming the Reparametrization Lemma}

Let  $\mu\in \mathcal{M}(M,T)$ and let $\gamma>0$. By Lemma \ref{thes} we choose $p_{\mu}\in \mathbb{N}\setminus\{0\}$ such that
\begin{eqnarray*}
\overline{\chi_1^+}(\mu)=\inf_{n\in \N}\frac{\int\log^+\|D_xT^n\|d\mu(x)}{n}&\geq &\frac{\int\log^+\|D_xT^{p_{\mu}}\|d\mu(x)}{p_{\mu}}-\frac{(r-1)\gamma}{4}. \label{ffirst}
\end{eqnarray*}

One can also assume $p_{\mu}$ large enough such that $H\left(p_{\mu}\frac{(r-1)\gamma}{4}\right)R(T)<\frac{(r-1)\gamma}{4}$ and $\frac{A+2\gamma}{p_{\mu}}<\frac{\gamma}{4}$. We will prove the statement of the Main Proposition with the entropy structure $\frac{1}{p_{\mu}}\mathcal{H}_{T^{p_\mu}}^{New}|_{\mathcal{M}(M,T)}$ (see Lemma \ref{struun}).

By continuity of $\mu\mapsto \int\log^+\|D_yT^{p_{\mu}}\|d\mu(y)$ and by upper semicontinuity of $\overline{\chi_1^+}$ one can choose the parameter $\delta_{\mu}>0$ such that for all ergodic measures $\nu$ with $\dist(\nu,\mu)<\delta_{\mu}$ we have :

 $$\left| \int \log^+\|D_{y}T^{p_{\mu}}\|d\nu(y)-\int\log^+\|D_yT^{p_{\mu}}\|d\mu(y) \right| < \frac{(r-1)\gamma}{4},$$

 $$\overline{\chi_1^+}(\mu)>\chi_1^+(\nu)-\frac{(r-1)\gamma}{4}.$$

We fix some ergodic measure $\nu$ with $\dist(\mu,\nu)<\delta_{\mu}$ with at most one nonnegative Lyapunov exponent. One can assume $h(\nu)>\frac{3\gamma}{4}$ because $\frac{1}{p_{\mu}}h^{New}_{T^{p_{\mu}}}(M|\nu,\epsilon_{\mu})\leq h(\nu)$ and the right member of Inequality (\ref{ggeq}),  $\frac{1}{r-1}\left(\overline{\chi_1^+}(\mu)-\chi_1^+(\nu)\right)+\gamma$,  is larger than $\frac{3\gamma}{4}$. In particular we have by Ruelle-Margulis inequality  $\chi_1^+(\nu)\geq h(\nu)>0$.
The measure $\nu$ need not be ergodic under $T^{p_{\mu}}$ but has at most $p_{\mu}$ ergodic components $\tilde{\nu}$ which all satisfy $\chi_i(\tilde{\nu},T^{p_{\mu}})=p_{\mu}\chi_i(\nu)$ for all $i=1,...,d$ and thus $l_{\tilde{\nu}}=l_{\nu}=1$.

We fix such an ergodic component $\tilde{\nu}$. For any $0<\alpha<1$ we let  $F_\alpha$ be the Borel  set  of $\tilde{\nu}$ measure larger than $\alpha$ as in Section \ref{volume}.

One can also assume by the ergodic theorem  that $\left(\lambda_n^+(x,T^{p_{\mu}})\right)_{n\in \N}$ are converging uniformly in $x\in F_{\alpha}$ to $ \int \log^+\|D_yT^{p_{\mu}}\|d\tilde{\nu}(y)$ and that $F_{\alpha}\subset G_{\alpha}$ with  $G_{\alpha}$ as in Section \ref{volume} for $T^{p_{\mu}}$ and $\tilde{\nu}$.\\

Let $\sigma:[0,1]\rightarrow M$ be a $1$-disk  associated to $F_{\alpha}$ and $\tilde{\nu}$ satisfying Inequality (\ref{vvolume}), i.e.
with $\epsilon$ and $D$ as in (\ref{vvolume})
  \begin{eqnarray*}
  H_{T^{p_{\mu}}}(n,\delta|x,F_{\alpha},\epsilon)&\leq & De^{\gamma n}V_{x,T^{p_{\mu}}}^{n,\epsilon}(\sigma|\chi_1(\tilde{\nu}),\gamma,C)
 \end{eqnarray*}

 We apply now the Reparametrization Lemma to $T^{p_{\mu}}$, $\sigma$  and to a given point $x\in F_{\alpha}$ : there exist $\epsilon>\epsilon_{\mu}>0$  depending only on $\|T^{p_{\mu}}\|_s$, $s=1,...,[r],r$, a real number $B$ depending only on $\chi_1(\tilde{\nu}),\gamma,C$ and $\|\sigma\|_s$, $\|T^{p_{\mu}}\|_s$, $s=1,...[r],r$, and families $(\mathcal{F}_n)_{n}$ of affine  maps from $[0,1]$ to $[0,1]$ satisfying the properties (i)-(ii)-(iii) of the Reparametrization Lemma. By (i) the volume of $T^{n-1}\circ \sigma \circ\psi$ is less than or equal to $1$ for all $\psi\in \mathcal{F}_n$. Therefore we have with (ii)

$$V_{x,T^{p_{\mu}}}^{n,\epsilon}(\sigma|\chi_1(\tilde{\nu}),\gamma,C)\leq \sharp \mathcal{F}_n$$

and  it follows from (iii) that :

\begin{multline*} h_{T^{p_{\mu}}}(M|F_{\alpha},\epsilon_{\mu})
\leq   \lim_{n\rightarrow +\infty}\sup_{x\in
F_{\alpha}}
\frac{1}{r-1}\left(1+H([\lambda^+_n(x,T^{p_{\mu}})-\chi_1(\tilde{\nu})]+3)\right)\left(\lambda^+_n(x,T^{p_{\mu}})-\chi_1(\tilde{\nu})\right)\\+ A+2\gamma
\end{multline*}

According to the definition of $F_{\alpha}$  we have :

\begin{eqnarray*}
\lim_{n\rightarrow +\infty}\inf_{x\in F_{\alpha}}\lambda^+_n(x,T^{p_{\mu}})-\chi_1(\tilde{\nu})&= & \lim_{n\rightarrow +\infty}\sup_{x\in F_{\alpha}}\lambda^+_n(x,T^{p_{\mu}})-\chi_1(\tilde{\nu})\\
& = & \int\log^+\|D_yT^{p_{\mu}}\|d\tilde{\nu}(y)-\chi_1(\tilde{\nu})\geq 0
\end{eqnarray*}

Now we distinguish cases :

\begin{itemize}
\item either $ \int\log^+\|D_yT^{p_{\mu}}\|d\tilde{\nu}(y)<\chi_1(\tilde{\nu})+p_{\mu}\frac{(r-1)\gamma}{4}$, then  the term \\
$\lim_{n\rightarrow +\infty}\sup_{x\in
F_{\alpha}} H([\lambda^+_n(x,T^{p_{\mu}})-\chi_1(\tilde{\nu})]+3)\left(\lambda^+_n(x,T^{p_{\mu}})-\chi_1(\tilde{\nu})\right)$
is bounded from above by $\log(2)\left(\int\log^+\|D_yT^{p_{\mu}}\|d\tilde{\nu}(y)-\chi_1(\tilde{\nu}) \right)$ which is less than  $\log(2)p_{\mu}\frac{(r-1)\gamma}{4}$.

\item or $\int\log^+\|D_yT^{p_{\mu}}\|d\tilde{\nu}(y)\geq \chi_1(\tilde{\nu})+p_{\mu}\frac{(r-1)\gamma}{4}$, then we have $\lim_{n\rightarrow +\infty}\inf_{x\in F_{\alpha}}\lambda^+_n(x,T^{p_{\mu}})-\chi_1(\tilde{\nu})\geq p_{\mu}\frac{(r-1)\gamma}{4}$. But we choose $p_{\mu}$ large enough so that $R(T)H\left(p_{\mu}\frac{(r-1)\gamma}{4}\right)<\frac{(r-1)\gamma}{4}$. It follows that :
\begin{eqnarray*}
&\lim_{n\rightarrow +\infty} \sup_{x\in F_{\alpha}}H([\lambda^+_n(x,T^{p_{\mu}})-\chi_1(\tilde{\nu})]+3)\left(\lambda^+_n(x,T^{p_{\mu}})-\chi_1(\tilde{\nu})\right)  \\
\leq & H\left(p_{\mu}\frac{(r-1)\gamma}{4}\right)p_{\mu}R(T) \\
\leq & p_{\mu}\frac{(r-1)\gamma}{4}
\end{eqnarray*}
\end{itemize}

We get finally in both cases  :

\begin{eqnarray*} h_{T^{p_{\mu}}}(M|F_{\alpha},\epsilon_{\mu})\leq  \frac{1}{r-1}\left(\int\log^+\|D_yT^{p_{\mu}}\|d\tilde{\nu}(y)- \chi_1(\tilde{\nu})\right)+p_{\mu}\frac{\gamma}{4}+A+2\gamma
\end{eqnarray*}

Then  by  letting $\alpha$ go to $1$ we obtain  since $\frac{A+2\gamma}{p_{\mu}}<\frac{\gamma}{4}$ :

$$\frac{1}{p_{\mu}}h_{T^{p_{\mu}}}^{New}(M|\tilde{\nu},\epsilon_{\mu})\leq
\frac{1}{(r-1)p_{\mu}}\left(\int\log^+\|D_yT^{p_{\mu}}\|d\tilde{\nu}(y)- \chi_1(\tilde{\nu})\right)+\frac{\gamma}{2} $$

The above inequality holds for all ergodic components $\tilde{\nu}$ of $\nu$ and then also for $\nu$ by harmonicity :
\begin{eqnarray}\label{nouvel}
\frac{1}{p_{\mu}}h_{T^{p_{\mu}}}^{New}(M|\nu,\epsilon_{\mu})&\leq &
\frac{1}{(r-1)p_{\mu}}\left(\int\log^+\|D_yT^{p_{\mu}}\|d\nu(y)- p_{\mu}\chi_1(\nu)\right)+\frac{\gamma}{2}
\end{eqnarray}

 Now we deduce from the choice of $\delta_{\mu}$ that :

$$\int \log^+\|D_yT^{p_{\mu}}\|d\nu(y)- p_{\mu}\chi_1(\nu)  \leq \int\log^+\|D_yT^{p_{\mu}}\|d\mu(y)- p_{\mu}\chi_1(\nu)+\frac{(r-1)\gamma}{4}$$

 and then by the choice of $p_{\mu}$ we get :

\begin{eqnarray*}
  \frac{1}{p_{\mu}}\left(\int\log^+\|D_yT^{p_{\mu}}\|d\nu(y)-p_{\mu}\chi_1(\nu)\right)& \leq & \overline{\chi^+_1}(\mu)-\chi_1^+(\nu)+\frac{(r-1)\gamma}{2}
\end{eqnarray*}

Together with Inequality (\ref{nouvel}) we conclude that
\begin{eqnarray*}
 \frac{1}{p_{\mu}}h_{T^{p_{\mu}}}^{New}(M|\nu,\epsilon_{\mu}) &\leq &\frac{1}{r-1}\left(\overline{\chi^+_1}(\mu)-\chi_1^+(\nu)\right)+\gamma.
\end{eqnarray*}
\end{demof}

The end of this paper deals with the proof of  the Reparametrization Lemma. We first state the key ingredients. The first one is a combinatorial  argument which allows us to work with subset of Bowen balls of length $n$ where the defect of multiplicativity of the norm $\frac{\|DT\|\|DT^k\circ \sigma\|}{\|DT^{k+1}\circ \sigma\|}$ is fixed for $k=0,...,n-1$. Then we recall Yomdin reparametrization lemma and a Landau-Kolmogorov inequality which  will be used to get $\mathcal{C}^r$ estimates of the Newhouse local entropy. Finally we explain how bound the local volume of a $\mathcal{C}^1$ curve by assuming that the derivative of this curve oscillates little.

\section{Technical lemmas}\label{feu}
This section is devoted to some useful technical lemmas for the proof of the Reparametrization Lemma presented in the last section.

\subsection{Combinatorial Lemma}

We  first begin with a usual combinatorial lemma which was already used by T.Downarowicz and A.Maass in \cite{Dow} and by the author in \cite{burinv}.

\begin{defi}
Let $S\in \N$ and $n\in\N$. We say that a sequence of $n$ positive integers $\mathcal{K}_n:=(k_1,...,k_n)$ admits the value $S$ if $\frac{1}{n}\sum_{i=1}^nk_i\leq S$.
\end{defi}

The number of  sequences of $n$ positive integers admitting the value $S$ is exactly the binomial coefficient ${nS \choose n}$. When $S$ is large enough this term is exponentially small in $nS$. More precisely we have the following lemma (recall that $H$ denotes the function defined from $[1,+\infty[$ by $H(t)=-\frac{1}{t}\log(\frac{1}{t})-(1-\frac{1}{t})\log(1-\frac{1}{t})$).

\begin{lem}\label{combi}
The logarithm of the number of sequences of $n$ positive integers admitting the value $S$ is at most
$nSH(S)+1$.
\end{lem}

We refer to  \cite{flup} (Lemma 16.19) for a proof. In the two next subsections we present the tools which allows us to get $\mathcal{C}^r$ estimates of the symbolic extension entropy function in contrast with \cite{burinv} where the author only deals with $\mathcal{C}^2$ maps.

\subsection{Estimates à la Yomdin}\label{landau}

 We recall now the heart of Yomdin theory which estimates the "local differential complexity" of intermediate smooth maps. The proof is based on a powerful Semi-algebraic Lemma due to M.Gromov \cite{Gr} (see also \cite{BurY} for a complete proof of Gromov's statement).

\begin{lem}\label{coeur}\cite{Gr}
Let $k$ and $d$ be positive integers and let $s>0$ be a positive real number. For any positive real number $a>0$ and for any $\mathcal{C}^s$ map $g:[0,1]^k\rightarrow \R^d$ with $\|g\|_s\leq a$,  there exists a family of real analytic maps $\mathcal{F}$ from $[0,1]^k$ to itself such that :

\begin{itemize}
 \item $\bigcup_{\phi\in\mathcal{F}}\phi ([0,1]^k)\supset \left\{x\in [0,1]^k, \ \|g(x)\|\leq a\right\}$,
 \item $\|\phi\|_t\leq 1$ for all $t=1,2,...,[s]+1$,
 \item $\|g\circ\phi\|_t\leq \frac{a}{12e}$ for $t=\min(1,s),2...,[s],s$,
 \item $\sharp \mathcal{F}\leq C$ with a universal constant $C$ depending only on $k,d$ and $s$.
\end{itemize}
\end{lem}

In fact this is a functional version of Lemma 3.4 of \cite{Gr} but the proofs are the same. The constant $\frac{1}{12e}$ may be replaced by any other universal constant. We choose this one for convenience of computation in the final proof of the Reparametrization Lemma.\\

The Semi-algebraic Lemma \cite{Gr} claims that the statement holds for polynomial functions without any condition on the derivatives, but the constant $C$ may depend on the degree. Then to prove the general case one approximates the given  intermediate smooth function by its Lagrangian polynomial. The hypothesis on the highest order derivative makes the  approximation good enough to conclude the proof. \\

To bound the local volume growth of $k$-disks  Y.Yomdin uses Lemma \ref{coeur} in a dynamical context. More precisely, given a $\mathcal{C}^r$ map $T:M\rightarrow M$ on a compact manifold $M$ and a $k$-disk $\sigma$ with $\|T^n\circ \sigma\|_s\leq 1$, he applies it with $a=1$, $s=r$ and $g=T^n\circ \sigma$.  Then it follows from (iii) that the $k$-volume of $T^n\circ\sigma$ is universally bounded from above. To ensure the condition on the $s$-norm we have to subdivide a general disk in $C_{ste}\|DT\|^{\frac{nk}{r}}$ subdisks so that the exponential rate of the local volume growth of a $k$-disk is bounded by $k\frac{R(T)}{r}$  replacing  if necessary $T$ by an iteration of $T$.\\

Here we use Lemma  \ref{coeur} only to reparametrize areas where the oscillation of the derivative of  $(T^n\circ \sigma)'$ is small compared to its size. Thus we will apply Lemma \ref{coeur} with $s=r-1$ to the derivative $(T^n\circ \sigma)'$. That's why we get the factor $\frac{1}{r-1}$ (and not $\frac{1}{r}$ as above) in the estimates of the Newhouse local entropy given by the Main Proposition. The constant $a$ will be chosen to be the norm of the derivative of $T^n\circ\sigma$ at typical points. Then in the case of $1$-disks an easy  geometrical argument presented in Subsection \ref{curve} allows us to bound the local volume growth.

\subsection{A Landau-Kolmogorov type Inequality}\label{landau}
Yomdin's reparametrization maps are semi-algebraic and it seems difficult to control from below the size of their derivative. Therefore  we just use affine maps in the proof of the Reparametrization Lemma. The key ingredient  to control the derivatives  is then the following classical Landau-Kolmogorov type inequality due to L.Neder \cite{neder}.

\begin{lem}\label{lan}\cite{neder}
Let $g:[0,1]\rightarrow \R^d$ be a $\mathcal{C}^{s}$ map with $s>0$. Then, there exists a universal  constant
$C$ depending only on $s$ and $d$ such that

$$\forall k=0,1,...,[s], \ \|g\|_{k}\leq C\left(\|g\|_{0}+\|g\|_{s}\right)$$
\end{lem}

For the sake of completeness we give a short proof of this result, that we borrow from  \cite{dev} (Theorem 5.6). In this reference the result is stated for integers $2\leq s$ but it can be easily extended to the general case $1<s\in \R$ as follows.

\begin{demo}
Without loss of generality we may assume  $d=1$. Let $x\in [0,\frac{1}{2}]$. From Taylor's formula, we have for all $0\leq u\leq \frac{1}{2}$ :

\begin{eqnarray*}
g(x+u)&=&g(x)+ug'(x)+...+\frac{u^{\left\lceil s-2\right\rceil}}{\left\lceil s-2\right\rceil !}g^{(\left\lceil s-2\right\rceil)}(x)+
\int_{0}^u\frac{(u-t)^{\left\lceil s-2\right\rceil}}{\left\lceil s-2\right\rceil !}f^{\left\lceil s-1\right\rceil}(x+t)dt \nonumber\\
&= & g(x)+ug'(x)+...+\frac{u^{\left\lceil s-2\right\rceil}}{\left\lceil s-2\right\rceil !}g^{(\left\lceil s-2\right\rceil)}(x)+\frac{u^{\left\lceil s-1\right\rceil}}{\left\lceil s-1\right\rceil !}g^{(\left\lceil s-1\right\rceil)}(x)\nonumber \\
& & \ \ \ \  \ \  \  \ \  +
\int_{0}^u\frac{(u-t)^{\left\lceil s-2\right\rceil}}{\left\lceil s-2\right\rceil !}\left(g^{\left\lceil s-1\right\rceil}(x+t)-g^{\left\lceil s-1\right\rceil}(x)\right)dt
\end{eqnarray*}

As $g^{(\left\lceil s-1\right\rceil)}$ is $s-\left\lceil s-1\right\rceil$ Hölder the integral remainder term  $$R(x,u)=\int_{0}^u\frac{(u-t)^{\left\lceil s-2\right\rceil}}{\left\lceil s-2\right\rceil !}\left(g^{\left\lceil s-1\right\rceil}(x+t)-g^{\left\lceil s-1\right\rceil}(x)\right)dt$$ can be bounded as follows

\begin{eqnarray}\label{labn}
\left|R(x,u)\right| & \leq &  \int_0^u\frac{(u-t)^{\left\lceil s-2\right\rceil}t^{s-\left\lceil s-1\right\rceil}}{\left\lceil s-2\right\rceil !} \|g\|_sdt\nonumber\\
&\leq &  \frac{u^s\int_0^1(1-x){\left\lceil s-2\right\rceil}x^{s-\left\lceil s-1\right\rceil}dx}{\left\lceil s-2\right\rceil !} \|g\|_sdt\nonumber \\
&\leq & \frac{u^{s}B(s+1-\left\lceil s-1\right\rceil,\left\lceil s-1\right\rceil)\|g\|_s}{\left\lceil s-1\right\rceil !}
\end{eqnarray}
where $B$ denotes the usual beta function. Then we choose arbitrarily real numbers $0<\lambda_1<...\lambda_{\left\lceil s-1\right\rceil}<1$ and we let $A$ be the Vandermonde matrix given by  $A:=\left(\lambda_i^j\right)_{1\leq i,j\leq \left\lceil s-1\right\rceil}$. We have

$$\left(g\left(x+\frac{\lambda_i}{2}\right)-g(x)-R\left(x,\frac{\lambda_i}{2}\right)\right)_i=A\left(\frac{g^{(i)}(x)}{i!2^i}\right)_i$$

As the determinant of $A$  is nonzero this system of equations can be solved and we have according to Inequality (\ref{labn}) for all $x\in [0,\frac{1}{2}]$

$$\frac{\left|g^{(i)}(x)\right|}{i!}\leq C\left(\|g\|_0+\|g\|_{s}\right)$$

where $C$ is a constant depending only on $s$ and on the real numbers $\lambda_1<...<\lambda_{\left\lceil s-1\right\rceil}$. We argue similarly for $x\in [\frac{1}{2},1]$ by considering $-\frac{1}{2}\leq u\leq 0$. This concludes the proof of the lemma.
\end{demo}

We will apply this lemma with $s=r-1$ to the derivative $(T^n\circ \sigma)'$ with $T^n:=T_{n}\circ ... \circ T_1$ for a sequence $(T_n)_n$ of $\mathcal{C}^r$ maps defined on the unit ball of $\R^d$ and for a smooth $1$-disk $\sigma:[0,1]\rightarrow \R^d$.

\subsection{Curves with small oscillations of the derivative}\label{curve}
In the present subsection we give an alternative strategy to bound the local volume of a $\mathcal{C}^1$ curve $\sigma:[0,1] \rightarrow \R^d$, that is the length of
$\sigma|_{\sigma^{-1}\left(B(0,1)\right)}$.  We assume now that we control the oscillation of the derivative of $\sigma$  (instead of the $r$-derivative in Yomdin's approach). More precisely we want the ratio $\frac{\|\sigma'(t)-\sigma'(s)\|}{\|\sigma\|_1}$ to be small  uniformly in $t,s\in [0,1]$. Then the following easy geometrical argument allows us to bound the local volume of $\sigma$. This idea was already exploited in \cite{burinv} under a slightly different form (Proposition 8 of \cite{burinv}). In the following the space $\R^d$ is always endowed with the usual euclidean norm $\|\|$ and $B(x,r)$ will denote the ball of radius $r$ centered at $x\in \R^d$.

\begin{lem}\label{oscille}
Let $\sigma:[0,1] \rightarrow \R^d$ be a $\mathcal{C}^1$ curve satisfying the following properties

\begin{enumerate}
\item $\sigma([0,1])\cap B(0,1)\neq \emptyset$,
\item $\|\sigma'(t)-\sigma'(s)\|\leq \frac{\|\sigma\|_1}{3}$ for all $t,s\in [0,1]$.\\
\end{enumerate}

Then, there exists $[a,b]\subset [0,1]$ such that

\begin{itemize}
\item $\sigma([0,1])\cap B(0,1)\subset \sigma([a,b])\subset B(0,\sqrt{d})$,
\item $(b-a)\|\sigma\|_1\leq \sqrt{3d}$.
\end{itemize}
\end{lem}

\begin{demo}
Let $w\in [0,1]$ with $\sigma(w)\in B(0,1)$. The hypothesis (2) on the derivative $\sigma'$ of $\sigma$ implies that $\sigma'(s)$ lies in the cone $C:=\{u\in \R^d, \ \left|\angle u,\sigma'(w)\right|\leq \frac{\pi}{6} \}$ for all $s\in [0,1]$. Indeed we have for all $s\in [0,1]$ :

\begin{eqnarray*}
\|\sigma'(s)\|\sin \left|\angle \sigma'(s),\sigma'(w)\right| & \leq & \| \sigma'(s)-\sigma'(w)\|\\
& \leq & \frac{\|\sigma\|_1}{3}\\
&\leq & \frac{\|\sigma'(s)\|}{2}
\end{eqnarray*}

where the last inequality follows, with $s_0\in [0,1]$ such that  $\|\sigma'(s_0)\|=\|\sigma\|_1$, from :

\begin{eqnarray*}
\|\sigma'(s)\|&\geq &\|\sigma'(s_0)\|-\|\sigma'(s)-\sigma'(s_0)\|\\
&\geq & \|\sigma\|_1-\frac{1}{3}\|\sigma\|_1=\frac{2}{3}\|\sigma\|_1
\end{eqnarray*}

\scalebox{0.45}{\includegraphics{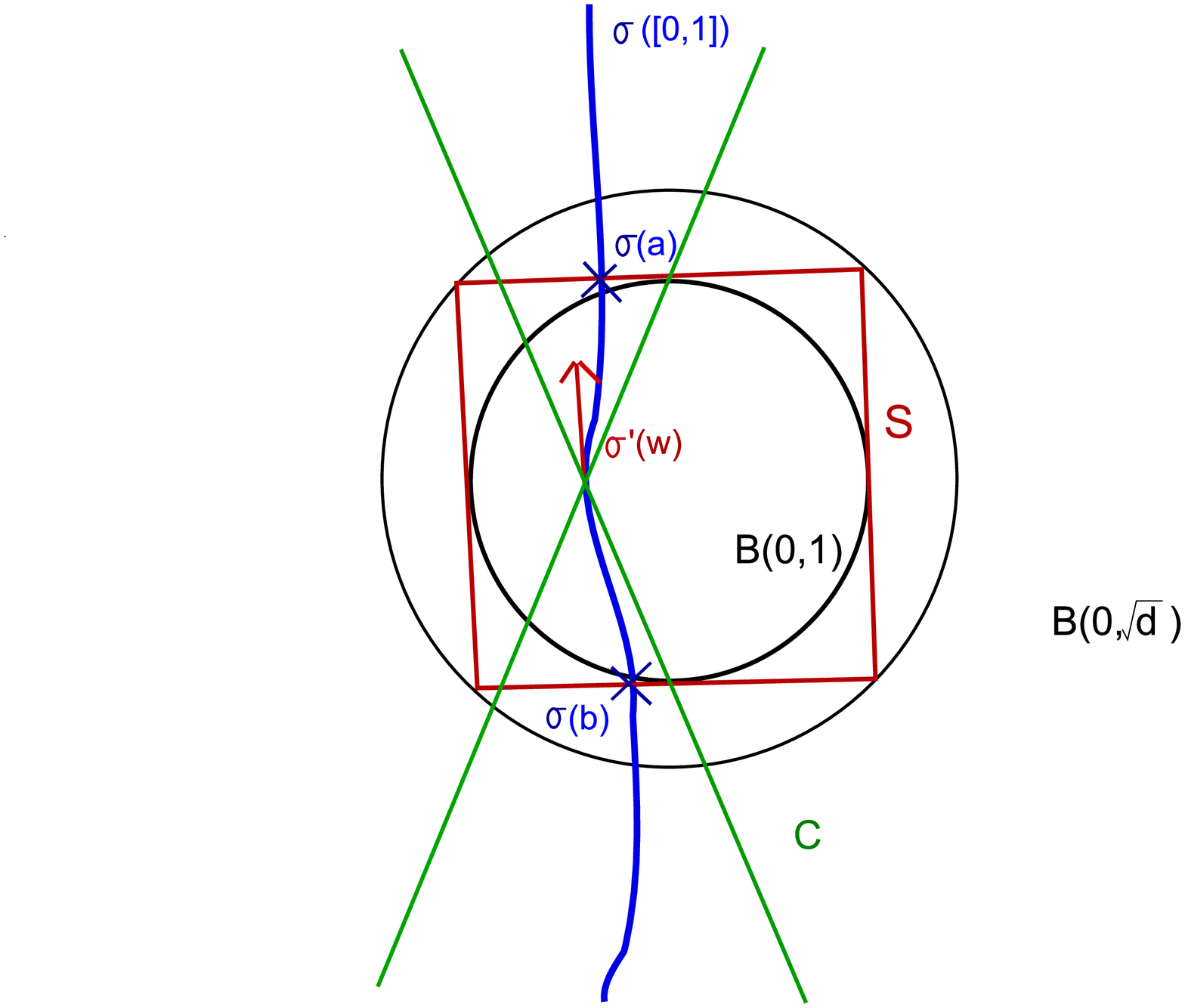}}

We consider the cube $S$ of size $\sqrt{d}$ containing the unit ball $B(0,1)$ whose faces are either orthogonal or parallel to $\sigma'(w)$. If $\sigma(t)$ does not belong to $S$ for some $t>w$ then $\sigma(s)$ stays  out of $S$ for $s>t$ (see the picture above). Similarly if $\sigma(t)$ is not in $S$ for some $t<w$ then so does $\sigma(s)$ for $s<t$. Therefore  if we set
$a=\sup\{s\leq w \ : \ \sigma(s)\in S \}$ and $b =\inf\{s\geq w \ : \ \sigma(s)\in S\}$ we have

$$\sigma([0,1])\cap B(0,1)\subset \sigma([a,b])\subset B(0,\sqrt{d})$$

Finally  we check the  second item :

\begin{eqnarray*}
(b-a)\|\sigma\|_1&\leq & \frac{3}{2}\int_a^b\|\sigma'(u)\|du \\
 &\leq & \sqrt{3}\int_a^b \frac{\sigma'(u).\sigma'(w)}{\|\sigma'(w)\|}du\\
 &\leq &\sqrt{3}\left\| \int_{a}^b\sigma'(u)du\right\|\\
 &\leq & \sqrt{3}\left\| \sigma(a)-\sigma(b)\right\| \leq \sqrt{3d}
\end{eqnarray*}

\end{demo}

\section{Reparametrization of $(n,\epsilon)$ Bowen balls : proof of the Reparametrization Lemma}

Let $M$ be a compact manifold of dimension $d$ and let $T:M\rightarrow M$ be a $\mathcal{C}^r$ map with $r>1$.
As in Yomdin theory, we consider the local dynamic at one point. We fix a Riemannian structure $\|\|$ on  $M$ and we denote by $R_{inj}$ the radius of injectivity and by $exp:TM(R_{inj})\rightarrow M$ the exponential map, where $TM(r):=\{(x,u), \ u\in T_xM, \ \|u\|_x<r\}$. We fix $R<R'<R_{inj}$ such that $T(B(x,R))\subset B(Tx,R')$ for all $x\in M$.\\

 Let $x\in M$ and $n\in \N$. We consider the map $T_n^x:T_{T^{n-1}x}M(R)\rightarrow T_{T^nx}M(R')$  defined by $T_n^x=exp^{-1}_{T^nx}\circ T\circ exp_{T^{n-1}x}$. For all $\epsilon<\frac{R}{\sqrt{d}}$, we put $T_{n,\epsilon}^x=\epsilon^{-1}T_n^x(\epsilon .):B(0,\sqrt{d})\rightarrow T_{T^nx}M\simeq \R^d$ and  $\mathcal{T}_{\epsilon}^x:=(T_{n,\epsilon}^x)_{n\in \N}$. Observe that for all small universal constant $A$ there exists $\epsilon>0$ depending only on  $\max_{s=\min(2,r),3,...,[r],r}\|D^sT\|$ such that $\|D^sT_{n,\epsilon}^x\|\leq A^{-1}$ for all $s=\min(2,r),2,...,[r], r$ and $n\in \N$.\\


 From now on we consider a general sequence $\mathcal{T}:=(T_n)_{n\in\N}$ of $\mathcal{C}^r$ maps with $r>1$ from $B(0,\sqrt{d})\subset \R^d$ to $\R^d$ with  $T_n(0)=0$ for all $n\in \N$. By convention we set $T_0=Id|_{B(0,\sqrt{d})}$. For each $n\in \N$ we  denote by $T^n$ the composition $T_n\circ...\circ T_0$ defined on $B(n,\sqrt{d})$ where $B(n,\rho)$ is the Bowen ball centered at $0$ of length $n$ and size $\rho>0$, i.e. $B(n,\rho):=\{y\in \R^d \ : \ \forall  k=0,..., n-1,\  \|T^ky\|<\rho\}$. \\

We extend the definition of $\mathcal{H}^n_{\mathcal{T}}$ to this framework by defining for all $1$-disks $\sigma$ and for all real numbers $\chi>0$, $\gamma>0$ and $C>1$ the set :
\begin{multline*}
\mathcal{H}^n_{\mathcal{T}}(\sigma,\chi,\gamma,C):=\Bigg\{t\in [0,1]\cap  \sigma^{-1}\left(B(n,\sqrt{d})\right) \ : \\ \forall 1\leq i\leq n, \ C^{-1}e^{(\chi-\gamma)i}\leq \|D_t\left(T^i\circ\sigma\right)\|\leq Ce^{(\chi+\gamma)i}\Bigg\}
\end{multline*}

We define now subsets of the Bowen ball $B(n+1,1)$ where the defect of multiplicativity of the norm of the composition $DT_{i+1}\circ D(T^i\circ \sigma)$ is prescribed at each step $1\leq i\leq n$.

\begin{defi} Let $\mathcal{T}:=(T_n)_{n\in\N}$ be a sequence of $\mathcal{C}^1$ maps from $B(0,\sqrt{d})\subset \R^d$ to $\R^d$ and let $\sigma:[0,1]\rightarrow \R^d$ be a $1$-disk of class $\mathcal{C}^1$. Let $\mathcal{K}_n:=(k_1,...,k_n)$  be a sequence of $n$ positive integers, we denote by $\mathcal{H}(\mathcal{K}_n)$ the subset of $\sigma^{-1}\left(B(n+1,1)\right)$ defined by :

\begin{multline}\mathcal{H}(\mathcal{K}_n):=\Bigg\{t\in [0,1]\cap \sigma^{-1}\left(B(n+1,1)\right)\ : \\ \forall 1\leq i\leq n, \ \left[ \log^+ \frac{\|D_t(T^i\circ \sigma)\|\max(1,\|D_{T^i\circ \sigma(t)}T_{i+1}\|)}{\|D_t(T^{i+1}\circ\sigma)\|}\right] +1=k_i \Bigg\}
\end{multline}

\end{defi}

 In the next lemma we estimate as in \cite{burinv} the numbers of such sets intersecting $\mathcal{H}^n_{\mathcal{T}}(\sigma,\chi,\gamma,C)$  by using the combinatorial argument of Lemma \ref{combi}.  We write  $\lambda^+_n:=\frac{1}{n}\sum_{i=0}^{n-1}\log^+\|D_{0}T_i\|$. When $\mathcal{T}$ is the sequence $\mathcal{T}_{\epsilon}^x$ for some $\epsilon>0$ and $x\in M$  we have $\lambda_n^+=\lambda_n^+(x,T)=\frac{1}{n}\sum_{i=0}^{n-1}\log^+\|D_{T^ix}T\|$. If we assume

 \begin{equation}\label{osd}
 \forall n\in \N \ \forall z,z'\in B(0,\sqrt{d}), \ \frac{\max(\|D_zT_n\|,1)}{\max(\|D_{z'}T_n\|,1)}\leq 2
 \end{equation}

then we have for all $y\in B(n,\sqrt{d})$

 \begin{eqnarray*}
 \left|\lambda^+_n-\frac{1}{n}\sum_{i=0}^{n-1}\log^+\|D_{y}T_i\|\right|&\leq & \log 2
 \end{eqnarray*}

The condition (\ref{osd}) is fulfilled as soon as the uniform norm of second derivative (or the $r-1$ Hölder norm when $r<2$) of each $T_n$ is less than $1/\sqrt{d}$. Under this condition we also notice that that for two points $t,s$ lying in $\mathcal{H}(\mathcal{K}_n)$ and satisfying $\frac{1}{2}\|D_t(T^n\circ\sigma)\|\leq \|D_s(T^n\circ\sigma)\|\leq 2\|D_t(T^n\circ\sigma)\|$  the derivative of $T^{n+1}\circ\sigma$ at $t$ and $s$ have also almost the same size. More precisely we have

 \begin{eqnarray}\label{lourd}
\frac{1}{4e}\|D_t(T^{n+1}\circ\sigma)\|\leq \|D_s(T^{n+1}\circ\sigma)\|\leq 4e\|D_t(T^{n+1}\circ\sigma)\|
 \end{eqnarray}

This last remark will be used in the proof of Proposition \ref{repint}. We state now the application of the combinatorial fact of Lemma \ref{combi}.

\begin{lem}\label{com}
Let $\chi>0$, $\frac{1}{3}>\gamma>0$ and $C>1$  and let $\mathcal{T}:=(T_n)_{n\in\N}$ be a sequence of $\mathcal{C}^r$ maps  from $B(0,\sqrt{d})$ to $\R^d$ such that $\frac{\max(\|D_zT_n\|,1)}{\max(\|D_{z'}T_n\|,1)}\leq 2$ for all $n\in\N$ and for all $z,z'\in B(0,\sqrt{d})$.\\

Then there exists an integer $N$ depending only on $C$ such that for all $n>N$ and for all $1$-disks $\sigma:[0,1]\rightarrow \R^d$ of class $\mathcal{C}^1$ the number of sequences $\mathcal{K}_{n-1}$ such that $\mathcal{H}(\mathcal{K}_{n-1})$ has a non empty intersection with $\mathcal{H}^n_{\mathcal{T}}(\sigma,\chi,\gamma,C)$ is bounded above  by

$$e^{3n-2}e^{(n-1)(\lambda^+_n-\chi)H([\lambda^+_n-\chi]+3)}$$

\end{lem}

\begin{demo}
Indeed if $t\in \mathcal{H}^n_{\mathcal{T}}(\sigma,\chi,\gamma,C)$, then $\sum_{i=0}^{n-1}\log^+\|D_{T^i\circ \sigma(t)}T_{i+1}\|-\log \|D_t(T^n\circ\sigma)\|\leq
n(\lambda^+_n+\log 2-\chi+\gamma)+\log C$. Thus the sequence $$\left(\left[ \log \frac{\|D_{t}(T^i\circ\sigma)\|\max(\|D_{T^i\circ \sigma(t)}T_{i+1}\|,1)}{\|D_t(T^{i+1}\circ \sigma)\|}\right]+1\right)_{i=1,...,n-1}$$ admits the value  $[\lambda^+_n-\chi]+3$ for $n>\frac{\log C}{1-\log 2-\frac{1}{3}}$.  We apply finally the combinatorial Lemma \ref{combi} to conclude the proof.
\end{demo}

 The Reparametrization Lemma follows now directly from Lemma \ref{com} and the following Proposition \ref{repint} applied to the sequences $\mathcal{T}=\mathcal{T}_{\epsilon}^x$ for all $x\in M$ and  for $\epsilon>0$ small enough.

\begin{prop}\label{repint}
Let $\mathcal{T}:=(T_n)_{n\in\N}$ be a
sequence of $\mathcal{C}^r$  maps with $r>1$  from $B(0,\sqrt{d})\subset \R^d$ to  $\R^d$ with $\sup_{n\in\N}\|T_n'\|<+\infty$ and $\sup_{\substack{n\in \mathbb{N},\\ 2\leq s\leq r}}\|T_n\|_s\leq A^{-1}$ where $A$ is a universal constant depending only on $r$ and $d$ which we specify later on in the proof. \\

 Then for all integers $n$, for all $1$-disks $\sigma$ of class $\mathcal{C}^r$ with $\max_{ s=1,...,[r], r}\|\sigma\|_s\leq 1$ and  for all sequences $\mathcal{K}_{n-1}=(k_1,...,k_{n-1})$ of $n-1$ positive integers there exists a family $\mathcal{G}_n$ of affine maps $\phi_n:[0,1]\rightarrow [0,1]$
satisfying the following properties :

\begin{enumerate}[(i)]
\item  $\forall \phi_n\in \mathcal{G}_n, \ \sigma\circ \phi_n([0,1])\subset B(n+1,\sqrt{d})$,
\item  $\forall \phi_n\in \mathcal{G}_n, \  \forall k=0,...,n, \ \forall  s=1,\min(2,r),3,...,[r],r,  \ \|T^{k}\circ \sigma\circ \phi_n\|_s\leq 1$,
\item $\forall \phi_n\in \mathcal{G}_n, \ \forall \min(1,r-1)\leq s\leq r-1, \   \|(T^n\circ \sigma)'\circ \phi_n\|_s\leq \frac{1}{3}\|(T^n\circ \sigma)'\circ \phi_n\|_0$,
\item  $ \mathcal{H}(\mathcal{K}_{n-1})\cap \sigma^{-1}\left(B(n+1,1)\right)\subset \bigcup_{\phi_n\in \mathcal{G}_n}\phi_n([0,1])$,
\item  $\log \sharp \mathcal{G}_n \leq B+An+\frac{1}{r-1}\sum_{i=1}^{n-1}k_i$,\\
 with $B$ depending only on $\sigma$.
\end{enumerate}
\end{prop}

When $d=1$ it follows from the property $\|(T^n\circ \sigma)'\circ \phi_n\|_{\min(1,r-1)}\leq \frac{1}{3}\|(T^n\circ \sigma)'\circ \phi_n\|_0$ that the derivative of $T^n\circ \sigma$ does not vanish on the image of each $\phi_n$, which lies thus in a monotone branch of $T^n\circ\sigma$. Therefore our result recovers the estimates on the number of monotone branches intersecting  $\mathcal{H}(\mathcal{K}_{n})$ obtained in \cite{Dow}.\\

The main difference with the reparametrization result presented in \cite{burinv} is the affine property of the reparametrization charts. Moreover, as these charts are in this paper one  dimensional this makes much easier the proof because the changes of charts $\phi_{n+1}\circ \phi_{n}^{-1}$ are just affine maps of $[0,1]$. Observe also that we do not assume any hyperbolicity condition on the reparametrized subset of the Bowen ball contrarily to Proposition 12 of \cite{burinv}. \\

The conditions on the reparametrization maps are stronger than in Yomdin theory where only (ii) is required. However we only consider here one dimensional disks - Yomdin's approach applies in any dimension - and we do not reparametrized the whole Bowen ball but only a subset with a fixed growth of the derivative.

\begin{demo}
We argue by induction on $n$. The initial step is easily checked.  We assume  the existence of the family $\mathcal{G}_n$ and we build $\mathcal{G}_{n+1}$. Let $\phi_n\in \mathcal{G}_n$. We cover the unit interval into $[e^{\frac{k_n}{r-1}}]+1$ subintervals of size $e^{-\frac{k_n}{r-1}}$. We reparametrize them from $[0,1]$ by affine contractions $\theta_{n+1}$. We let $\psi_{n+1}:=\phi_n\circ\theta_{n+1}$.
From now on we focus on   maps $\psi_{n+1}$ whose image intersects the set
$$\mathcal{H}(\mathcal{K}_n)\cap \sigma^{-1}\left(B(n+2,1)\right)$$

Fix such a map $\psi_{n+1}$ and choose $w\in [0,1]$ such that $\psi_{n+1}(w)$ belongs to the previous set.  Let us first prove that

$$\|(T^{n+1}\circ \sigma)'\circ \psi_{n+1}\|_{r-1}\leq 2\|(T^{n+1}\circ \sigma)'\circ \psi_{n+1}(w)\|$$

To simplify the notations we write $T'_{n+1}$ for the differential $DT_{n+1}$ of $T_{n+1}$. Then, according to the chain rule derivative we have

$$(T^{n+1}\circ \sigma)'\circ \psi_{n+1}=T_{n+1}'\circ \left(T^n\circ\sigma\circ\psi_{n+1}\right)
\times (T^{n}\circ\sigma)'\circ \psi_{n+1}$$

By the induction hypothesis (iii) the $s$-norm of $(T^{n}\circ\sigma)'\circ \psi_{n+1}$ for  $\min(1,r-1)\leq s\leq r-1$ satisfies

$$\left\|(T^{n}\circ\sigma)'\circ \psi_{n+1}\right\|_s\leq \frac{1}{3}\|\theta_{n+1}\|_1^s \|(T^{n}\circ \sigma)'\circ \phi_{n}\|_0$$

We consider now the first term $T_{n+1}'\circ \left(T^n\circ\sigma\circ\psi_{n+1}\right)$ of the product. We first recall that \begin{itemize}
\item for an integer $\beta$ the $\beta$-derivative of a product or a composition of smooth functions is an universal polynomial in the $\alpha$-derivatives of each term with integers $\alpha=0,...,\beta$,
\item for a positive real number $0<\beta\leq 1$  the $\beta$-Hölder norm $\|f\times g\|_{\beta}$ of a product $fg$  is less than or equal to $\|f\|_{0}\|g\|_{\beta}+\|g\|_{0}\|f\|_{\beta}$ while the $\beta$-Hölder norm $\|f\circ g\|_{\beta}$ of a composition $f\circ g$ is less than or equal to
    $\|f\|_{\beta}\|g\|_{1}^{\beta}$.
\end{itemize}

Therefore,  by the induction hypothesis (ii), we have for all $\min(1,r-1)\leq s \leq r-1$
$$\left\|T_{n+1}'\circ \left(T^n\circ\sigma\circ \psi_{n+1}\right)\right\|_s\leq C\|\theta_{n+1}\|_1^s\sup_{\min(2,r)\leq k\leq r}\|T_{n+1}\|_k$$   where $C$ is a constant depending only on $r$ and $d$.  Finally we get by replacing this constant by another one which we denote again by $C$ :

\begin{eqnarray*}
\|(T^{n+1}\circ \sigma)'\circ \psi_{n+1}\|_{r-1}&\leq &  \|\theta_{n+1}\|_1^{r-1}\| \|(T^{n}\circ \sigma)'\circ\phi_{n}\|_0\Bigg(\frac{1}{3}\|T'_{n+1}\circ \left(T^n\circ \sigma\circ\psi_{n+1}\right)\|_0\\
& & \ \  \ \ \ \ \  \ \ \ \ \ \ \ \ \ \ \ \ \ \ \ \ \ \ \ \ \  \ \ \ \ \ \ \ \ \ +C\sup_{ \min(2,r)\leq k\leq r}\|T_{n+1}\|_k\Bigg)\\
&\leq & e^{-k_n}\|(T^{n}\circ \sigma)'\circ\phi_{n}\|_0\max(\|T'_{n+1}\|_0,1)\left(\frac{1}{3}+C\sup_{ \min(2,r)\leq k\leq r}\|T_{n+1}\|_k\right)
\end{eqnarray*}

But it follows from the induction hypothesis (iii) with $s=\min(1,r-1)$ that $$\|(T^{n}\circ \sigma)' \circ \psi_{n+1}(w)\|\geq  \frac{2}{3}\|(T^{n}\circ \sigma)' \circ \phi_{n}\|_0$$
and we can choose the constant $A$ large enough to ensure firstly
$$C\sup_{\substack{p\in \mathbb{N},\\ \min(2,r)\leq k\leq r}}\|T_p\|_{k}\leq CA^{-1} \leq  \frac{1}{3}$$ and secondly for all $x\in B(0,\sqrt{d})$
$$\max(\|T'_{n+1}\|_{0},1)\leq 2\max(1,\|T'_{n+1}(x)\|)$$

Therefore we have

\begin{eqnarray*}
\|(T^{n+1}\circ \sigma)' \circ \psi_{n+1}\|_{r-1}& \leq & \frac{2}{3}e^{-k_n}\|(T^{n}\circ \sigma)' \circ \phi_{n}\|_0 \max(\|T'_{n+1}\|_0,1)\\
&\leq& e^{-k_n}\|(T^{n}\circ \sigma)' \circ \psi_{n+1}(w)\| \max(\|T'_{n+1}\|_0,1)\\
&\leq &  2e^{-k_n}\|(T^{n}\circ \sigma)'\circ\psi_{n+1}(w)\|  \max(\|T'_{n+1}\left(T^n\circ\sigma\circ \psi_{n+1}(w)\right)\|,1)
\end{eqnarray*}

and since  $\psi_{n+1}(w)$ belongs to $\mathcal{H}(\mathcal{K}_n)$ we obtain

\begin{eqnarray} \label{dess}
\|(T^{n+1}\circ \sigma)' \circ \psi_{n+1}\|_{r-1} & \leq & 2\|(T^{n+1}\circ\sigma)' \circ \psi_{n+1}(w)\|
\end{eqnarray}

  Now we apply Lemma \ref{coeur}  to the $\mathcal{C}^{r-1}$ map $(T^{n+1}\circ\sigma)' \circ \psi_{n+1}$ with $a=4e\|(T^{n+1}\circ\sigma)' \circ \psi_{n+1}(w)\|$. We let $\mathcal{F}_{\psi_{n+1}}$ be the associated family of  reparametrization maps.
For any $\xi_{n+1}\in \mathcal{F}_{\psi_{n+1}}$
we have in particular
$$\|(T^{n+1}\circ\sigma)' \circ \psi_{n+1}\circ\xi_{n+1}\|_{\min(1,r-1)}\leq \frac{1}{3}\|(T^{n+1}\circ\sigma)' \circ \psi_{n+1}(w)\|$$

From now on we only consider the maps $\xi_{n+1}\in \mathcal{F}_{\psi_{n+1}}$ such that the image of $\psi_{n+1}\circ \xi_{n+1}$ meets $\mathcal{H}(\mathcal{K}_n)$.
 Since $\|(T^n\circ \sigma)'\circ \phi_{n}\|_{\min(1,r-1)}\leq \frac{1}{3}\|(T^n\circ \sigma)'\circ \phi_{n}\|_0$ which as seen earlier implies $\|(T^n\circ \sigma)'\circ \phi_{n}(t)\|\leq 2 \|(T^n\circ \sigma)'\circ \phi_{n}(s)\|$ for all $t,s\in [0,1]$ we have by the inequalities (\ref{lourd})

\begin{eqnarray*}
\frac{1}{8e}\|(T^{n+1}\circ\sigma)' \circ \psi_{n+1}(w)\|&\leq\|(T^{n+1}\circ\sigma)' \circ \psi_{n+1}\circ\xi_{n+1}\|_{0}&\leq  8e\|(T^{n+1}\circ\sigma)' \circ \psi_{n+1}(w)\|
\end{eqnarray*}

Let $\eta_{n+1}:[0,1]\rightarrow \xi_{n+1}([0,1])\subset [0,1]$ be an  affine reparametrization of $\xi_{n+1}([0,1])$. It follows from the above Landau-Kolmogorov Inequality of Lemma \ref{lan} that for all $\min(1,r-1)\leq s \leq r-1$ we have

\begin{eqnarray*}
\|(T^{n+1}\circ\sigma)' \circ \psi_{n+1}\circ\eta_{n+1}\|_s&\leq & C\big(\|(T^{n+1}\circ\sigma)' \circ \psi_{n+1}\circ\eta_{n+1}\|_0\\
& & \ \  \ \ \ \ \  \ \ \ \ \ \ \ \  \ \ \ \ \ \  \ \ \ \ \ \ \ \ \ +\|(T^{n+1}\circ\sigma)' \circ \psi_{n+1}\circ\eta_{n+1}\|_{r-1}\big)\\
&\leq & C\left(\|(T^{n+1}\circ\sigma)' \circ \psi_{n+1}\circ\xi_{n+1}\|_0+\|(T^{n+1}\circ\sigma)' \circ \psi_{n+1}\|_{r-1}\right)\\
&\leq & C\|(T^{n+1}\circ\sigma)' \circ \psi_{n+1}(w)\|\\
&\leq & C\|(T^{n+1}\circ\sigma)' \circ \psi_{n+1}\circ\eta_{n+1}\|_0
\end{eqnarray*}

where the universal constant $C$ may change at each step of the previous sequence of inequalities.

By dividing the unit interval into $[3C]+1$ subintervals of size $<1/3C$ and by reparametrizing them affinely one can assume  $\|(T^{n+1}\circ \sigma)'\circ\psi_{n+1}\circ\eta_{n+1}\|_s\leq \frac{1}{3}\|(T^{n+1}\circ \sigma)'\circ\psi_{n+1}\circ\eta_{n+1}\|_0$.

Now for any map $T^{n+1}\circ\sigma \circ \psi_{n+1}\circ\eta_{n+1}$ we let $[a_{\eta_{n+1}},b_{\eta_{n+1}}]$ be the subinterval of $[0,1]$ given by Lemma \ref{oscille}. We reparametrize $\psi_{n+1}\circ \eta_{n+1}|_{[a_{\eta_{n+1}},b_{\eta_{n+1}}]}$ from $[0,1]$ by an affine contraction to get new maps that we denote by $\phi_{n+1}$, i.e. $\phi_{n+1}(t):=
\psi_{n+1}\circ \eta_{n+1}\left(a_{\eta_{n+1}}+(b_{\eta_{n+1}}-a_{\eta_{n+1}})t\right)$ for all $t\in [0,1]$. By construction the family $\mathcal{G}_{n+1}$ of  affine maps $\phi_{n+1}$ satisfy  properties (i) and (iv). Moreover, since the map $\psi_{n+1}\circ\eta_{n+1}$ satisfies (iii) then so does $\phi_{n+1}$. Therefore, we only need to check (ii) for the family $\mathcal{G}_{n+1}$ at step $n+1$.

By Lemma \ref{oscille} we have  $\|T^{n+1}\circ \sigma\circ\phi_{n+1}\|_1\leq \sqrt{3d}$. It only remains  to check (ii) for the family $\mathcal{G}_{n+1}$ at step $n+1$ for the derivative of order $\min(2,r)\leq s\leq r$. By using successively  the affine property of $\phi_{n+1}$, property (iii) and  the above case $s=1$, we have for $\min(2,r)\leq s\leq r$ :

\begin{eqnarray*}
\|T^{n+1}\circ \sigma \circ \phi_{n+1}\|_s&=&\|(T^{n+1}\circ \sigma)' \circ \phi_{n+1}\|_{s-1}\times \|\phi_{n+1}\|_1\\
& \leq & \frac{1}{3}\|(T^{n+1}\circ \sigma)' \circ \phi_{n+1}\|_0\times \|\phi_{n+1}\|_1\\
&\leq & \frac{1}{3}\|T^{n+1}\circ \sigma \circ \phi_{n+1}\|_1\leq \sqrt{\frac{d}{3}}
\end{eqnarray*}
 As previously done we may suppose this last constant to be one, should we multiply the number of reparametrization maps by $\left[\sqrt{\frac{d}{3}}\right]+1$. This concludes the proof  of the lemma.
\end{demo}

\it{E-mail  address :}  \texttt{David.Burguet@cmla.ens-cachan.fr}

\end{document}